\documentclass[10pt,twoside]{siamart1116}

\usepackage[english]{babel}
\usepackage{graphicx,epstopdf,epsfig}
\usepackage{amsfonts,epsfig,fancyhdr,graphics, hyperref,amsmath,amssymb}
\usepackage{longtable}

\newtheorem{remark}[theorem]{Remark}
\newtheorem{example}[theorem]{Example}
\newtheorem{observation}[theorem]{Observation}
\newtheorem{question}[theorem]{Question}

\renewcommand{\theequation}{\thesection.\arabic{equation}}
\renewtheorem{theorem}{Theorem}[section]

\usepackage{etaremune, enumerate, float, verbatim}
\usepackage{algorithm}
\usepackage{algorithmic}
\usepackage{theoremref}
\usepackage{graphicx}
\usepackage{multirow} 
\usepackage{soul} 
\usepackage{url}
\usepackage[mathlines]{lineno}
\usepackage{dsfont} 
\usepackage{tikz, graphicx, subcaption, caption}
\usepackage[algo2e,ruled,vlined]{algorithm2e}
\usepackage{mathrsfs,amsmath}
\usepackage{color}
\definecolor{red}{rgb}{1,0,0}

\definecolor{blue}{rgb}{0,0,.7}

\definecolor{green}{rgb}{0,.6,0}

\definecolor{purp}{rgb}{.5,0,.5}

\numberwithin{figure}{section}
\numberwithin{table}{section}
\usepackage{theoremref}

\newcommand{\RKm}[1]{R(K_{#1},\vec K_{#1})  }
\newcommand{\SKm}[1]{S(K_{#1},\vec K_{#1})  }
\newcommand{\orth}{\mathcal{O}} 
\newcommand{\Q}{\mathcal{Q}} 

\newcommand{\R}{\mathbb{R}}

\newcommand{\sym}{\operatorname{sym}}
\newcommand{\sgn}{\operatorname{sgn}}
\newcommand{\znz}{\operatorname{znz}}
\newcommand{\diag}{\operatorname{diag}}
\newcommand{\rank}{\operatorname{rank}}
\newcommand{\pll}{{\mkern3mu\vphantom{\perp}\vrule depth 0pt\mkern4mu\vrule depth 0pt\mkern3mu}} 
\newcommand{\pert}{r}   

\newcommand{\x}{\times}

\newcommand{\bit}{\begin{itemize}}
\newcommand{\eit}{\end{itemize}}
\newcommand{\ben}{\begin{enumerate}}
\newcommand{\een}{\end{enumerate}}
\newcommand{\beq}{\begin{equation}}
\newcommand{\eeq}{\end{equation}}
\newcommand{\bea}{\begin{eqnarray*}}
\newcommand{\eea}{\end{eqnarray*}}
\newcommand{\bean}{\begin{eqnarray}}
\newcommand{\eean}{\end{eqnarray}}
\newcommand{\bpf}{\begin{proof}}
\newcommand{\epf}{\end{proof}\ms}
\newcommand{\bmt}{\begin{bmatrix}}
\newcommand{\emt}{\end{bmatrix}}
\newcommand{\ms}{\medskip}

\newcommand{\beqa}{\begin{array}}
\newcommand{\eeqa}{\end{array}}

\newcommand{\lf}{\left\lfloor}
\newcommand{\rf}{\right\rfloor}
\newcommand{\lp}{\left(}
\newcommand{\rp}{\right)}
\newcommand{\lb}{\left[}
\newcommand{\rb}{\right]}

\newcommand{\wh}{\widehat}

\newcommand{\bc}{{\bf c}}
\newcommand{\bv}{{\bf v}}
\newcommand{\bw}{{\bf w}}

\newcommand{\bx}{{\bf x}}
\newcommand{\by}{{\bf y}}
\newcommand{\br}{{\bf r}}
\newcommand{\bq}{{\bf q}}

\newcommand{\bzero}{{\bf 0}}
\newcommand*{\abs}[1]{\lvert #1\rvert}
\newcommand*{\squig}[1]{\widetilde{#1}}
\DeclareMathOperator*{\Span}{span}
\newcommand*{\norm}[1]{\lVert #1\rVert}
\newcommand*{\what}[1]{\widehat{#1}}

\DeclareMathOperator*{\supp}{supp}

\newcommand{\Names}{Zachary Brennan, Christopher Cox, Bryan A.~Curtis, Enrique Gomez-Leos, Kimberly P.~Hadaway, Leslie Hogben, Conor Thompson}
\newcommand{\Title}{
Orthogonal realizations of random sign patterns and other applications of the SIPP}

\author{Zachary Brennan\thanks{Department of Mathematics, Iowa State University (brennanz, cocox, bcurtis1, enriqueg, kph3, hogben, conorjt)@iastate.edu}\and Christopher Cox\footnotemark[2]\and Bryan A.~Curtis\footnotemark[2]\ \thanks{Corresponding author}\and Enrique Gomez-Leos\footnotemark[2] \and  Kimberly P.~Hadaway\footnotemark[2]\and Leslie Hogben\footnotemark[2]\ \thanks{American Institute of Mathematics, San Jose, CA 95112, USA (hogben@aimath.org).}\and Conor Thompson\footnotemark[2]}

 \title{\Title}

\begin{document}
\maketitle

\begin{abstract}
A \emph{sign pattern} is an array with entries in $\{+,-,0\}$.  
A matrix $Q$ is \emph{row orthogonal} if $QQ^T = I$. The Strong Inner Product Property (SIPP), introduced in [B.A.~Curtis and B.L.~Shader, Sign patterns of orthogonal matrices and the strong inner product property, \emph{Linear Algebra  Appl.} 592: 228--259, 2020], is an important tool when   determining whether a sign pattern allows row orthogonality because it guarantees there is a nearby matrix with the same property, allowing zero entries to be perturbed to nonzero entries, while preserving the sign of every nonzero entry. This paper uses the SIPP to initiate the study of conditions under which random sign patterns allow row orthogonality with high probability. Building on  prior work,  $5\x n$ nowhere zero sign patterns that minimally allow orthogonality are determined.  Conditions on zero entries  in a sign pattern are established   that guarantee any row orthogonal matrix with such a sign pattern has the SIPP.
\end{abstract}



%
%
\section{Introduction}

A \emph{sign pattern} is an array with entries coming from the set $\{+,-,0\}$. The entries of sign patterns 
encode qualitative properties of real matrices. Sign patterns were introduced in applications where the entries of the matrix may be known only approximately (or not at all), but the signs of the entries are known. A matrix $Q$ is \emph{row orthogonal} provided $QQ^T = I$. The problem of determining whether an $m\x n$ sign pattern allows row orthogonality has been studied for many years \cite{FM1964,JW98,CS2020,Curtis_20}. Recently the strong inner product property (SIPP) was introduced by Curtis and Shader in \cite{CS2020} to study sign patterns of row orthogonal matrices. This paper relies heavily on the SIPP to build on prior work (e.g., classifying small patterns that minimally allow orthogonality) and initiate the study of conditions under which random sign patterns allow row orthogonality with high probability.
 
Finding a certificate that a sign pattern allows row orthogonality is often difficult. By applying a variant of Gram-Schmidt orthogonalization to a nowhere zero nearly row orthogonal matrix we obtain conditions that guarantee the existence of a nearby row orthogonal matrix with the same sign pattern (see Section \ref{s:approx-orthog}).  
In Section \ref{s:SIPP}, we  apply the SIPP to develop new tools to show that a sign pattern allows row orthogonality and use these tools (and the results from Section \ref{s:approx-orthog}) to determine $5\x n$ nowhere zero sign patterns that minimally allow orthogonality.  We also establish conditions on zero entries in a sign pattern  that guarantee any row orthogonal matrix with such a sign pattern has the SIPP.
One of our main results, Theorem \ref{thm:main prob allows}, utilizes the SIPP to obtain a lower bound $h(m)$ such that the probability of a random $m\x n$ sign pattern allowing row orthogonality goes to $1$ as $m$ tends toward $\infty$ for $n\ge h(m)$ (here random means $+$ and $-$ are equally likely and the probability of $0$ is given).
The remainder of this introduction defines terminology and notation  (Section \ref{s:defs}) and lists known results we will use (Section \ref{s:resultscited}).  

\subsection{Definitions and notation}\label{s:defs} In the study of sign patterns, sometimes the distinction between zero and nonzero is more important than the sign.  A \emph{zero-nonzero pattern} or \emph{znz pattern} is an array with entries coming from the set $\{*,0\}$.  We use the term \emph{pattern} to mean a sign pattern or a zero-nonzero pattern. 
Given a real number $a$,
\[
\sgn(a) =
\begin{cases}
+ & \text{if } a>0\\
- & \text{if } a<0\\
0 & \text{if } a=0
\end{cases}
\qquad \text{and} \qquad
\znz(a) =
\begin{cases}
* & \text{if } a\not=0\\
0 & \text{if } a=0.
\end{cases}
\]
The sign pattern and zero-nonzero pattern of a matrix $A = [a_{ij}]$ are $\sgn(A) = [\sgn(a_{ij})]$ and $\znz(A) = [\znz(a_{ij})]$, respectively. The qualitative class of an $m\x n$ sign pattern $S$ is the set
\[
\Q(S) = \{A \in \R^{m\times n} : \sgn(A) = S\},
\]
and the qualitative class of an $m\times n$ znz pattern $Z = [z_{ij}]$ is the set
\[
\Q(Z) = \{A \in\R^{m\times n} : \znz(A) = Z\}.
\]
A matrix in the qualitative class $\Q(P)$ is called a \textit{realization} of the pattern $P$.  For a sign pattern $S$,  $C_S$ denotes the unique $(1,-1,0)$-matrix that is a realization of the sign pattern $S$. Similarly, $C_Z$ is the unique $(1,0)$-matrix that is a realization of the zero-nonzero pattern $Z$.  
A \textit{superpattern} of a sign pattern $S=[s_{ij}]$ is a sign pattern $R=[r_{ij}]$ of the same dimensions such that $r_{ij}=s_{ij}$ whenever $s_{ij}\in\{+,-\}$; if $s_{ij}=0$ then $r_{ij}\in\{+,-,0\}$.

A matrix with orthogonal rows is not necessarily row orthogonal; for us, the rows of a row orthogonal matrix have unit length. The set of $m \x n$ row orthogonal matrices is denoted by $\orth(m,n)$ and  we write $\orth(m)$ as shorthand for $\orth(m,m)$. Note that  every matrix $Q\in \orth(m)$  is orthogonal, i.e., $Q^TQ=QQ^T=I$.  The set of $m\x m$ real symmetric matrices is denoted by $\sym(m)$.

A \emph{zero matrix} $O\in\R^{m\x n}$ or \emph{zero vector} $\bzero\in\R^n$ has every entry equal to zero.  An $m\x n$ matrix or  pattern is \emph{wide} if $m\le n$.  A wide matrix has \emph{full rank} if its rank equals its number of rows, i.e., it has linearly independent rows. A row orthogonal matrix is necessarily wide. Let  $A\in\R^{m\times n}$ be a wide matrix. Then $A$ has the \textit{strong inner product property} (SIPP) provided $X=O$ is the only symmetric matrix satisfying $(XA)\circ A=O$ \cite{Curtis_20}. The strong inner product property is one of a number of strong properties of matrices that guarantee  there is a nearby matrix with the same property, allowing zero entries to be perturbed to nonzero entries, while preserving the sign of every nonzero entry  \cite[Part 2]{HLS22}.

An $m\x n$ sign pattern $S$ \emph{allows row orthogonality}  if there is a row orthogonal matrix $Q\in \Q(S)$ (equivalently, $\orth(m,n)\cap \Q(S)\ne\emptyset$). An $m\x n$ sign pattern $S$ \emph{allows o-SIPP} if there is a row orthogonal matrix $Q\in \Q(S)$ that has the SIPP. Since scaling a matrix with a positive constant does not change its pattern, no pattern requires row orthogonality. A sign pattern $S$ {\em requires o-SIPP} if  every $Q\in \Q(S)\cap\orth(m,n)$  has the SIPP and $ \Q(S)\cap\orth(m,n)\ne \emptyset$. Without the assumption that $ \Q(S)\cap\orth(m,n)\ne \emptyset$, the all zeros pattern would vacuously require o-SIPP.

An $m \times n$ rectangular sign pattern $S$ is row \emph{potentially pairwise-orthogonal} or \emph{row PPO}  if no row is a zero vector  and for each pair $(i,k)$ with  $1\le i<k\le m$, there are realizations of row $i$ and row $k$ that are  orthogonal. The term \emph{column PPO} is defined analogously. 
A pair of rows $i$ and $k$ in a $m \times n$ rectangular sign pattern $S=[s_{ij}]$ \emph{has a negative $4$-cycle} if there are two columns $j$ and $\ell$ such that $s_{ij}s_{kl}=+$ and $s_{ij}s_{kl}=-$, where multiplication on the set $\{0,+,-\}$ is defined in the obvious way that conforms to real arithmetic.
A pair of rows $i$ and $k$  in an $m\x n$  matrix or pattern $P=[p_{ij}]$   is \emph{combinatorially orthogonal} if $p_{ij}\ne 0$ implies $p_{kj}=0$ for every $j=1,\dots,n$. 
  
A \textit{signed permutation matrix} is a square $(1,-1,0)$-matrix with exactly one nonzero entry in each row and column. Matrices $A,B \in \mathbb{R}^{m \times n}$ are \textit{sign equivalent} if $A = P_1 B P_2$, where $P_1$ and $P_2$ are signed permutation matrices. Two sign patterns $S$ and $S'$ are sign equivalent if $C_S$ and $C_{S'}$ are sign equivalent.

For a vector $\bv\in \R^n$, the \emph{support} of $\bv$, denoted by $\supp(\bv)$, is the set of indices of nonzero entries of $\bv$.   Let $[n]=\{1,\dots,n\}$.

\subsection{Known results}\label{s:resultscited} 
In the remainder of this introduction we provide some known results about the SIPP that we will use.
The primary motivation for developing the SIPP is given by the next theorem of Curtis and Shader.\footnote{Theorem 4.5 in \cite{CS2020} actually says that every superpattern of $S$ allows row orthogonality, but the proof shows it allows o-SIPP.} We provide a slightly stronger result in Theorem \ref{thm: gen sup pat}.

\begin{theorem}\label{t:basic-SIPP}{\rm \cite{CS2020}}
If $Q\in\orth(m,n)$ has the SIPP and $\sgn(Q)=S$, then every superpattern of $S$ allows o-SIPP.
\end{theorem}

Theorem \ref{t:basic-SIPP} has many consequences. Here we list some {that} we use. A matrix with the SIPP or a sign pattern that allows the SIPP can be padded with additional zero columns and retain that property.

\begin{lemma}{\rm \cite{CS2020}}\label{lem: add O} 
Let $A\in\R^{m\times n}$ and $p>n$. Then $A$ has the SIPP if and only if the $m\times p$ matrix $B=\big[ A~|~O\big]$ has the SIPP.\end{lemma}

\begin{corollary}{\rm \cite{CS2020}}\label{thm: SIPP add zero block} 
If $Q\in\orth(m,n)$ has the SIPP and $\sgn(Q)=S$,  then $\bmt{S~\big\vert~ O}\emt$ allows o-SIPP. 
\end{corollary}

The next two results show that sign equivalence preserves having the SIPP, as does taking the transpose for (square) orthogonal matrices.

\begin{proposition}{\rm \cite{CS2020}}\label{prop: signequiv}
Let $A,B\in\R^{m\times n}$ be sign equivalent. Then $A$ has the SIPP if and only if $B$ has the SIPP.
\end{proposition}

\begin{proposition}{\rm \cite{CS2020}}\label{prop: equiv transpose}
Let $Q\in\orth(m)$. Then $Q$ has the SIPP if and only if $Q^T$ has the SIPP.
\end{proposition}

The previous results provide some sufficient conditions for a sign pattern to allow row orthogonality. The next result provides a way to show a sign pattern does not allow row orthogonality. 

\begin{theorem}\label{thm: big rank 1}{\rm \cite{JW98}}
Let $S$ be a nowhere zero sign pattern and let $R$ be an $r\times s$ submatrix of $S$. If $r+s \geq n+2$ and $\rank C_R = 1$, then $S$ does not allow row orthogonality.
\end{theorem}


%
%

\section{From approximate orthogonality to exact orthogonality}\label{s:approx-orthog}

In this section, we establish a result that gives conditions under which a collection of ``nearly'' orthogonal vectors necessarily implies the existence of a ``nearby'' collection of truly orthogonal vectors.
Such a result is similar in spirit to the effective implicit function theorems used by, e.g., Cohn, Kumar and Minton \cite{Cohn2016} to derive the existence of an exact code from an approximate one. However, instead of using an implicit function theorem, we simply rely on the Gram--Schmidt process.  Although the perturbations here are created by a different mechanism, we also point out that the idea of perturbing one solution to obtain another desired solution is a fundamental idea underlying strong properties.  First, we define the  function  $\pert_m(\epsilon)$ used to quantify the notion of ``nearby.''  

\begin{definition}\label{d:r_m} For an integer $m\ge 1$ and a real number $0\leq\epsilon<\frac 1 {m-1}$, 
\[
\pert_m(\epsilon)=\sqrt{{1+\epsilon\over (1-(m-2)\epsilon)(1-(m-1)\epsilon)}}-1
\]
where $\frac 1 0$ is interpreted as $\infty$, so $\pert_1(\epsilon)=0$ for all $\epsilon\geq 0$.
\end{definition}

For simplicity, we have defined the functions $\pert_m(\epsilon)$ in closed form.  In order to use these functions for the results of this section we need a recursive approach, which is given in  the next lemma.

\begin{lemma}\label{lem:rm-recursive}
Given $\pert_1(\epsilon)=0$,   $\pert_m(\epsilon)$ can be computed recursively for  all $m\geq 2$ and all $0\leq\epsilon<\frac 1 {m-1}$ by  
\[
\pert_m(\epsilon) = \sqrt{\frac{1+\epsilon}{ 1-\epsilon}}\biggl(\pert_{m-1}\biggl(\frac{\epsilon}{1-\epsilon}\biggr)+1\biggr)-1.
\]
\end{lemma}
\begin{proof} It is easy to verify the result for $m=2$. 
For $m\geq 3$,
\begin{align*}
\pert_m(\epsilon) &=  \sqrt{{1+\epsilon\over(1-(m-2)\epsilon)(1-(m-1)\epsilon)}}-1.\\
&= \sqrt{{1+\epsilon\over 1-\epsilon}} \sqrt{{(1-\epsilon)((1-\epsilon)+\epsilon)\over \bigl((1-\epsilon)-(m-3)\epsilon \bigr)\bigl((1-\epsilon)-(m-2)\epsilon\bigr)}}-1\\
&= \sqrt{{1+\epsilon\over 1-\epsilon}} \sqrt{{1+{\epsilon\over 1-\epsilon}\over \bigl(1-(m-3){\epsilon\over 1-\epsilon}\bigr)\bigl(1-(m-2){\epsilon\over 1-\epsilon}\bigr)}}-1\\
&=  \sqrt{{1+\epsilon\over 1-\epsilon}}\biggl(\pert_{m-1}\biggl({\epsilon\over 1-\epsilon}\biggr)+1\biggr)-1.
\end{align*}
\end{proof}

The next lemma provides the key step to go from approximately to exactly orthogonal.

\begin{lemma}\label{lem:project}
Let $m$ be a positive integer, let $0\leq\epsilon<{1\over m-1}$ and fix any  inner product space $(\Omega,\langle\cdot,\cdot\rangle)$. Additionally, let $\norm\cdot$ be any norm on $\Omega$ (possibly unrelated to $\langle\cdot,\cdot\rangle$). If $\bx_1,\dots,\bx_m\in\Omega$ satisfy
\begin{enumerate}
\item $\langle \bx_i,\bx_i\rangle=1$ for all $i\in[m]$, and
\item $\abs{\langle\bx_i,\bx_j\rangle}\leq\epsilon$ for all $i\neq j\in[m]$,
\end{enumerate}
then there exists $\squig\bx_1,\dots,\squig\bx_m\in\Span\{\bx_1,\dots,\bx_m\}$ satisfying
\begin{enumerate}
\item $\{\squig\bx_1,\dots,\squig\bx_m\}$ is orthonormal with respect to $\langle\cdot,\cdot\rangle$, and
\item $\norm{\bx_i-\squig\bx_i}\leq\pert_m(\epsilon)\norm{\bx_i}$ for all $i\in[m]$.
\end{enumerate}
\end{lemma}
\begin{proof}
We prove the result by induction on $m$. The case $m=1$ is immediate by taking $\squig \bx_1 = \bx_1$. Let $m\geq 2$ be a positive integer and suppose the statement holds for $m-1$.

Without loss of generality, $\norm{\bx_m}\leq \norm{\bx_i}$ for all $i \in\{1,\dots,m\}$.
For $i\in\{1,\ldots,m-1\}$, let $\bx_i^\pll=\langle \bx_i,\bx_m\rangle \bx_m$ and $\bx_i^\perp=\bx_i - \bx_i^\pll$; then $\langle\bx_i^\perp,\bx_m\rangle = 0$.
Since $\abs{\langle \bx_i,\bx_m\rangle}\leq\epsilon$ and $\langle\bx_m,\bx_m\rangle=1$, we know that $ 0\le \langle\bx_i^\pll,\bx_i^\pll\rangle\leq\epsilon^2$.
Therefore, the Pythagorean Theorem allows us to conclude that
\begin{equation}\label{eqn:2norm}
\langle\bx_i,\bx_i\rangle=\langle\bx_i^\pll,\bx_i^\pll\rangle+\langle\bx_i^\perp,\bx_i^\perp\rangle\quad\implies\quad 1 \ge   \langle\bx_i^\perp,\bx_i^\perp\rangle\geq 1-\epsilon^2.
\end{equation}
Since $\epsilon < 1$, we know that $\bx_i^\perp$ is non-zero. Let $\what \bx_i^\perp$ denote the unit vector in the direction of $\bx_i^\perp$, i.e., $\what\bx_i^\perp=\frac 1 {\sqrt{\langle\bx_i^\perp,\bx_i^\perp\rangle}}\bx_i^\perp$.

Since $\bx_i^\pll=\langle \bx_i,\bx_m\rangle \bx_m$, we see that  $\norm{\bx_i^\pll}\leq\epsilon\norm{\bx_m}\leq\epsilon\norm{\bx_i}$.
Then the triangle inequality applied to $\bx_i^\perp =\bx_i+(-\bx_i^\pll)$ implies that $\norm{\bx_i^\perp}\leq (1+\epsilon)\norm{\bx_i}$.
Together with \eqref{eqn:2norm}, we have
\begin{equation}\label{eqn:infnorm}
\norm{\what \bx_i^\perp}\leq\frac{1+\epsilon}{\sqrt{1-\epsilon^2}}\norm{\bx_i}=\sqrt{{1+\epsilon\over 1-\epsilon}}\norm{\bx_i}.
\end{equation}
In particular,
\bean
\norm{\bx_i-\what \bx_i^\perp} &\leq& \norm{\bx_i^\pll}+\norm{\bx_i^\perp-\what \bx_i^\perp}= \norm{\bx_i^\pll}+\biggl({1\over\sqrt{\langle\bx_i^\perp,\bx_i^\perp\rangle}}-1\biggr)\norm{\bx_i^\perp} \nonumber \\
&\leq& \epsilon\norm{\bx_i}+\biggl({1\over\sqrt{1-\epsilon^2}}-1\biggr)(1+\epsilon)\norm{\bx_i}  =\biggl(\sqrt{\frac{1+\epsilon}{ 1-\epsilon}}-1\biggr)\norm{\bx_i},\label{eqn:distance}
\eean
where the second inequality follows by combining \eqref{eqn:2norm} and $\norm{\bx_i^\perp}\leq (1+\epsilon)\norm{\bx_i}$.

Now, for any other $j\in\{1,\ldots,m-1\}$ with $j\neq i$, we have
\bea
\langle \bx_i^\perp, \bx_j^\perp\rangle &=&\langle \bx_i,\bx_j\rangle- \langle \bx_i,\bx_m\rangle\langle \bx_m,\bx_j\rangle\\
\implies \abs{\langle\bx_i^\perp,\bx_j^\perp\rangle} &\leq& \epsilon+\epsilon^2\\
\implies \abs{\langle \what \bx_i^\perp,\what \bx_j^\perp\rangle} &\leq& {\epsilon+\epsilon^2\over\sqrt{\langle\bx_i^\perp,\bx_i^\perp\rangle\langle\bx_j^\perp,\bx_j^\perp\rangle}}\leq\frac{\epsilon+\epsilon^2}{ 1-\epsilon^2} = \frac{\epsilon}{ 1-\epsilon}.
\eea

Therefore, since $\what\bx_1^\perp,\dots,\what\bx_{m-1}^\perp$ are unit vectors by construction, we may apply the induction hypothesis to find an orthonormal set $\{\squig \bx_1,\dots,\squig \bx_{m-1}\}\subseteq\Span\{\what\bx_1^\perp,\dots,\what\bx_{m-1}^\perp\}$ such that
\[
\norm{\what \bx_i^\perp-\squig \bx_i}\leq \pert_{m-1}\biggl(\frac{\epsilon}{ 1-\epsilon}\biggr)\norm{\what \bx_i^\perp}
\]
for each $i\in\{1,\ldots,m-1\}$.
By invoking additionally \eqref{eqn:infnorm} and \eqref{eqn:distance}, we bound
\bea
\norm{\bx_i-\squig \bx_i} &\leq\norm{\bx_i-\what \bx_i^\perp}+\norm{\what \bx_i^\perp-\squig \bx_i} \leq \biggl(\sqrt{{1+\epsilon\over 1-\epsilon}}-1\biggr)\norm{\bx_i}+\pert_{m-1}\biggl({\epsilon\over 1-\epsilon}\biggr)\norm{\what\bx_i^\perp}\\
&\leq\biggl(\sqrt{{1+\epsilon\over 1-\epsilon}}-1+\sqrt{{1+\epsilon\over 1-\epsilon}}\pert_{m-1}\biggl({\epsilon\over 1-\epsilon}\biggr)\biggr)\norm{\bx_i}=\pert_m(\epsilon)\norm{\bx_i}
\eea
for all $i\in\{1,\ldots,m-1\}$ where the last equality follows from Lemma \ref{lem:rm-recursive}.

Finally, let $\squig \bx_m = \bx_m$.  Then $\squig \bx_1,\dots,\squig \bx_m$ satisfy the claim, because $\squig \bx_1,\dots,\squig \bx_{m-1}\in\break\Span\{\what \bx_1^\perp,\dots,\what \bx_{m-1}^\perp\} \subseteq\Span\{\bx_1,\dots,\bx_m\}$ by construction, and the former subspace is orthogonal to $\bx_m=\squig\bx_m$.
\end{proof}

Observe that the process used to create the vectors $\squig \bx_i$ is a reordering of the modified Gram-Schmit process. We stated Lemma~\ref{lem:project} very generally in the hopes that other researchers will find it useful; for our uses, we specialize to the standard Euclidean inner product and the $\infty$-norm to attain a result related to row orthogonal realizations.  

We apply Lemma \ref{lem:project} to obtain 
Theorem \ref{c:close-enough}, which will be used in Section \ref{ss:5xn} to characterize $5\times n$ nowhere zero sign patterns that minimally allow orthogonality. 
Essentially this result says that for any matrix that is close to being row orthogonal,  there exists a nearby matrix  that is row orthogonal and has the same sign pattern.

\begin{definition}{\rm 
For a non-zero vector $\bx=[x_1,\dots,x_n]^T\in\R^n$, define
\[
\delta(\bx)={\min_{i\in[n]}\abs{x_i}\over\max_{j\in[n]}\abs{x_j}}={\min_{i\in[n]}\abs{x_i}\over\norm{\bx}_\infty}.
\]
}
\end{definition} 

\begin{theorem}\label{c:close-enough}\label{thm:approx}
Let $\bx_1,\dots,\bx_m\in\R^n$ be any non-zero vectors and let
$
\epsilon = \max_{i\neq j} \big\lvert \big\langle{ \bx_i \over \norm{\bx_i}_2}, {\bx_j \over \norm{\bx_j}_2} \big\rangle \big\rvert,
$
where $\langle\cdot,\cdot\rangle$ is the standard Euclidean inner product. If
\begin{enumerate}
\item $\epsilon<{1\over m-1}$, and
\item $\pert_m(\epsilon)<\min_{i\in[m]}\delta(\bx_i)$,
\end{enumerate}
then there exists an orthogonal set $\{\squig\bx_1,\dots,\squig\bx_m\}\subseteq\R^n$ satisfying $\sgn(\bx_i)=\sgn(\squig\bx_i)$ for all $i\in[m]$.
\end{theorem}

\begin{proof}
We apply Lemma~\ref{lem:project} to the vectors $ \frac1{\|\bx_i\|_2} \bx_i$, specializing to the Euclidean inner product and the $\infty$-norm, to locate an orthonormal set $\{\squig\bx_1',\dots,\squig\bx_m'\}\subseteq\R^n$ such that
\[
\biggl\lVert {\bx_i\over\norm{\bx_i}_2}-\squig\bx_i'\biggr\rVert_\infty\leq \pert_m(\epsilon)\biggl\lVert{\bx_i\over\norm{\bx_i}_2}\biggr\rVert_\infty\quad\implies\quad\norm{\bx_i-\norm{\bx_i}_2\cdot\squig\bx_i'}_\infty\leq\pert_m(\epsilon)\norm{\bx_i}_\infty
\]
for all $i\in[m]$. In particular, setting $\squig\bx_i=\norm{\bx_i}_2\cdot\squig\bx_i'$ for each $i\in[m]$, we know that $\{\squig\bx_1,\dots,\squig\bx_m\}$ is an orthogonal set and that
\[
\norm{\bx_i-\squig\bx_i}_\infty\leq r_m(\epsilon)\norm{\bx_i}_\infty<\delta(\bx_i)\norm{\bx_i}_\infty=\min_{j\in[n]}\abs{(\bx_i)_j}.
\]
Since $\abs{x-y}<\abs x\implies\sgn(x)=\sgn(y)$ for any $x,y\in\R$, we conclude that $\sgn(\bx_i)=\sgn(\squig\bx_i)$ for all $i\in[m]$.
\end{proof}

One particularly useful feature of Theorem~\ref{thm:approx} is that it can be used to present reasonable certificates of the existence of row orthogonal realizations; in fact, it implies that integer-valued certificates can always be found.
We illustrate this in the following example (which will be used in Section \ref{ss:5xn}).

\begin{example}\label{ex:5by6allows}
{\rm 
Consider the sign-pattern
\[
S=\begin{bmatrix}
 - & - & - & + & + & + \\
 + & + & - & + & + & + \\
 + & + & + & - & - & + \\
 + & + & + & + & + & - \\
 + & + & + & + & + & +
\end{bmatrix}.
\]
Explicitly writing down a row orthogonal realization of $S$ would be difficult since this requires exact arithmetic. Despite this, it is not too difficult for a computer to find realizations of $S$ that are row orthogonal up to floating-point error. For example, the following matrix is such a realization for $S$:
\[
A' =
\begin{bmatrix}
 -0.0743294 & -0.668965 & -0.222988 & 0.371647 & 0.0743294 & 0.594635 \\
 0.118415 & 0.59468 & -0.665382 & 0.0360018 & 0.195624 & 0.387344 \\
 0.511869 & 0.0620542 & 0.206774 & -0.259068 & -0.646593 & 0.454076 \\
 0.665978 & 0.0389929 & 0.0396191 & 0.681063 & 0.0657504 & -0.291912 \\
 0.02319 & 0.264691 & 0.660817 & 0.0611589 & 0.541388 & 0.442585 \\
\end{bmatrix}.
\]
Of course, $A'$ is not actually a row orthogonal matrix and so it does not directly demonstrate that $S$ has a row orthogonal realization; however, $A'$ does satisfy the hypotheses of Theorem~\ref{thm:approx}. In fact, by scaling and truncating $A'$ appropriately, we find the following integer-valued matrix which satisfies the hypotheses of Theorem~\ref{thm:approx} as well:
\[
A = \begin{bmatrix}
-8 & -74 & -25 & 41 & 8 & 65 \\
13 & 65 & -73 & 4 & 22 & 43 \\
56 & 7 & 23 & -28 & -71 & 50 \\
73 & 4 & 4 & 75 & 7 & -32 \\
3 & 29 & 73 & 7 & 60 & 49
\end{bmatrix}.
\]
Here and in similar examples, we  use $\delta$ to denote $\min_{i\in[m]}\delta(\br_i)$ where the $\br_i$ are the rows of the matrix. To apply Theorem~\ref{thm:approx} to the matrix $A$, observe that the value $\delta = \frac3{73}>.004$ is obtained from row $5$ of $A$ and the value $\epsilon={71\over\sqrt{146335965}}<0.006$ is obtained from rows $1$ and $4$. Since $\pert_5$ is increasing on its domain, $\pert_5(\epsilon) < \pert_5(0.006) < 0.03 < \delta$. We may therefore apply Theorem~\ref{thm:approx} to conclude that there exists a row orthogonal matrix with the same sign-pattern as $A$.

We will use these same basic ideas  in Section~\ref{ss:5xn} to write down reasonable certificates for the existence of row orthogonal realizations for other sign patterns.
}
\end{example}

%
%

\section{Results on the SIPP}\label{s:SIPP}

In this section we present results related to the SIPP and orthogonality. Section \ref{sub:tools}  contains some useful tools for studying matrices that have the SIPP. Of particular interest is Theorem \ref{thm: gen sup pat} which extends Theorem \ref{t:basic-SIPP}. In Section \ref{sub:requires} we investigate sign patterns that require o-SIPP. Section \ref{ss:5xn} utilizes the SIPP to provide a complete characterization of nowhere zero $m\times n$ sign patterns that  minimally allow orthogonality for $m\leq 5$. 

\subsection{Tools}\label{sub:tools}

Recall that an $m\times n$ matrix $A$ has the SIPP provided $O$ is the only symmetric matrix $X$ satisfying $(XA)\circ A = O$. It is often much easier to  construct a matrix with orthogonal rows as opposed to a row orthogonal matrix. The next lemma allows us to study row orthogonal matrices with the SIPP without first normalizing the rows.

\begin{lemma}\label{lem: scale rows}
Suppose $Q$ is an $m\x n$ full rank matrix with orthogonal rows that  has the SIPP and  $D$ is any $m\x m$ diagonal matrix with every diagonal entry nonzero.  Then  $DQ$ has the SIPP. Furthermore, $D$ can be chosen so that $DQ$ is row orthogonal.
\end{lemma}
\begin{proof}
Let $X \in \sym(m)$ such that $(XDQ) \circ (DQ)=O$.
By algebraic manipulation,
\[
O=(XDQ) \circ (DQ) = (DXD)Q \circ Q.
\]
Since $DXD \in \sym(m)$ and $Q$ has the SIPP,  it follows that $DXD=O$, which implies $X=O$.  Thus $DQ$ has the SIPP. Let $\br_i^T$ denote the $i$th row of $Q$. Define  $D=\diag\lp\frac 1{\|\br_1\|},\dots,\frac 1{\|\br_m\|}\rp$ and  $\wh{Q} = DQ$, so $\wh{Q}\in\orth(m,n)$.\,  
\end{proof}

The next three lemmas showcase additional hypotheses on $A$ that imply various entries in $X$ are $0$.

\begin{lemma}\label{lem:zero-rows-of-X}
Let  $A\in\R^{m\times n}$ be a wide matrix with full rank and let $X\in\R^{m\times n}$. Suppose that every entry of row $k$ of $A$ is nonzero.  If $(XA) \circ A = O$, then every entry of row $k$ of $X$ is zero.
\end{lemma}
\begin{proof}
Suppose $(XA)\circ A = O$. Let $\br_1^T, \dots, \br_m^T$ denote the rows of $X$. Then $(XA) \circ A = O$ implies that  $\br_k^T A = \bzero^T$. Since $A$ has full rank there exists a matrix $B$ such that $AB = I$ and so $\bzero^T = \br_k^T AB =  \br_k^T$.
\end{proof}

\begin{lemma}\label{lem:identity-circ-Q}
Suppose $Q\in\orth(m,n)$,  $X\in\sym(m)$ and $(XQ) \circ Q = O$. Then $I\circ X = O$. 
\end{lemma}
\bpf
Let $Y=XQ$ and write $Y=[y_{ij}]$, $X=[x_{ij}]$,  and $Q=[q_{ij}]$. Since $Q$ is row orthogonal,  $X=YQ^{T}$. The condition that $(XQ) \circ Q = O$ implies that  $y_{ij} =0$ if $q_{ij} \not =0$. Therefore,
\[
x_{ii} = (YQ^{T})_{ii}  = \sum_{j=1}^{m} y_{ij} q_{ij}  =0.
\]
In other words, $I \circ X = O$.
\epf

\begin{lemma}\label{lem:2nz-col}
Let $Q = [q_{ij}]\in\orth(m,n)$ and  $X = [x_{ij}]\in\sym(m)$ satisfy the equation $(XQ) \circ Q = O$. Suppose that the only two nonzero entries in column $j$ of $Q$ are $q_{ij}$ and $q_{kj}$. Then $x_{ik}=x_{ki}=0$.
\end{lemma}
\bpf
Since $(XQ) \circ Q = O$, 
\[
0 = ((XQ)\circ Q)_{ij} = \lp\sum_{\ell=1}^m x_{i\ell} q_{\ell j}\rp q_{ij} = (x_{ii}q_{ij}+x_{ik}q_{kj})q_{ij}.
\]
Since $q_{ij}\ne 0$, $x_{ii}q_{ij}+x_{ik}q_{kj}=0$.  By Lemma \ref{lem:identity-circ-Q}, $x_{ii} = 0$, so $x_{ik}q_{kj}=0$. Since $q_{kj}\ne 0$, $x_{ik}=0$. 
\epf

The next result extends one direction of \cite[Proposition 3.9]{CS2020}.
\begin{lemma}\label{lem:(1,1)-block-has-SIPP}
Suppose $A$ is a wide matrix  partitioned as a $2\x 2$ block matrix 
$A = \left[\begin{array}{c|c}
A_{1} & A_{2} \\
\hline
A_{3} & A_{4} \\
\end{array}\right]$ 
with $A_{3}, A_{4}$ both nowhere zero (or vacuous). If $A_{1}$ has the SIPP and $A$ is full rank, then $A$ has the SIPP.
\end{lemma}
\begin{proof}
Let 
$X = \left[\begin{array}{c|c}
X_{1} & X_{2} \\ 
\hline 
X_{2}^{T} & X_{4} 
\end{array}\right]$ 
be a symmetric matrix such that $(XA)\circ A = O$. By Lemma \ref{lem:zero-rows-of-X}, 
$\left[\begin{array}{c|c}
X_{2}^{T} & X_{4}
\end{array}\right] =O$. 
Therefore, 
$X = \left[ \begin{array}{c|c}
X_{1} & O \\ 
\hline 
O & O 
\end{array}\right]$. 
The equation $(XA) \circ A =O$ implies that $(X_{1} A_{1}) \circ A_{1} =O$. Since $X_{1}$ is symmetric and $A_{1}$ has the SIPP,  $X_{1} =O$ and $A$ has the SIPP.
\end{proof}

Manifold theory, and in particular having manifolds interesect transversally, plays a fundamental role in strong properties, including the SIPP; see \cite{HLS22} for more information. Smooth manifolds $\mathcal M$ and $\mathcal N$,  both in $\R^d$, intersect \emph{transversally} at a point $\bx$ if $\bx \in \mathcal M \cap \mathcal N$ and the intersection of the normal spaces of $\mathcal M$ at $\bx$ and of $\mathcal N$ at $\bx$ contains only the zero vector.  As the next result shows, a matrix $Q\in\orth(m,n)$ having the SIPP is equivalent to $\Q(\sgn(Q))$ and $\orth(m,n)$ intersecting transversally at $Q$.

\begin{theorem}\label{t:basic-SIPP2}{\rm \cite[Theorem 4.5]{CS2020}}
Let $Q\in\orth(m,n)$ have sign pattern $S$. The manifolds $Q(S)$ and $\orth(m, n)$ intersect transversally at $Q$ if and only if $Q$ has the SIPP. If $Q$ has the SIPP, then every superpattern of $S$ allows o-SIPP.
\end{theorem}

Theorem \ref{thm: gen sup pat} improves the previous result by allowing us to control the relative magnitudes of the formerly zero entries in $Q$ when applying the SIPP.  This requires the following theorem of van der Holst, Lov\'{a}sz, and Schrijver.

\begin{theorem}{\rm \cite{Holst1999}} \label{perturb_manifolds}
Let $\mathcal{M}_1(s)$ and $\mathcal{M}_2(t)$ be smooth families of manifolds in $\R^d$, and assume that $\mathcal{M}_1(0)$ and $\mathcal{M}_2(0)$ intersect transverally at $\by_0$. Then there is a neighborhood $W\subseteq\R^2$ of the origin and a continuous function $f:W\to\R^d$ such that $f(0,0) = \by_0$ and for each $\epsilon = (\epsilon_1,\epsilon_2)\in W$, $\mathcal{M}_1(\epsilon_1)$ and $\mathcal{M}_2(\epsilon_2)$ intersect transversally at $f(\epsilon)$.
\end{theorem}

Note that the statement of Theorem \ref{perturb_manifolds} applies to a more general setting than we require. For our purposes, one of the smooth families of manifolds is replaced with a manifold. In such a setting we may think of $f$ as a continuous function of one variable  from an interval about the origin to $\R^d$.

\begin{theorem}\label{thm: gen sup pat}
Let $Q\in\orth(m,n)$ have sign pattern $S$. If $Q$ has the SIPP, then for all $A\in\R^{m\times n}$ with $A\circ Q = O$,  for every  $\epsilon$ sufficiently small  there is a  matrix $M_\epsilon\in\mathcal{Q}(S)$ such that $M_\epsilon + \epsilon A\in\orth(m,n)$. Moreover, $M_\epsilon + \epsilon A$ has the SIPP.  
\end{theorem}
\begin{proof}
Suppose that $Q = [q_{ij}]$ has the SIPP and let $A = [a_{ij}]\in\R^{m\times n}$ satisfy $A\circ Q = O$. Define the smooth family of manifolds $\mathcal{M}_A(t)$ by
\[
\mathcal{M}_A(t) = \{ B = [b_{ij}]\in \R^{m\times n}: \sgn(b_{ij}) = \sgn(q_{ij}) \text{ if } q_{ij}\not=0, \text{ and } b_{ij} = a_{ij}t \text{ if } q_{ij} = 0\}
\]
for $t\in(-1,1)$. Since $Q$ has the SIPP, $\orth(m,n)$ and $\mathcal{M}_A(0) = \Q(S)$ intersect transversally at $Q$ by Theorem \ref{t:basic-SIPP2}. By Theorem \ref{perturb_manifolds}, there exists a continuous function $f:(-1,1)\to\R^{m\times n}$ such that $f(0) = Q$ and the manifolds $\mathcal{M}_A(\epsilon)$ and $\orth(m,n)$ intersect transversally at $f(\epsilon)$ for each $\epsilon>0$ sufficiently small. Since $f$ is continuous we may choose $\epsilon$ small enough so that   $M_\epsilon := f(\epsilon)\circ  C_Z \in \Q(S)$, where $Z$ is the zero-nonzero pattern of $Q$ and $C_Z$ is the unique $(1,0)$-matrix in $\Q(Z)$. To complete the proof, observe that $f(\epsilon) =  M_\epsilon + \epsilon A$. Moreover, $f(\epsilon)$ has the SIPP  by  Theorem \ref{t:basic-SIPP2} since $\mathcal{M}_A(\epsilon)$ and $\orth(m,n)$ intersect transversally at $f(\epsilon)$.
\end{proof}

We apply the previous theorem to prove the next result.
\begin{proposition}
Let
\[
S =
\left[\begin{array}{c|c}
S_1 & O \\
\hline
S_3 & S_4
\end{array}\right]
\]
be a sign pattern that allows row orthogonality, and let $S_4'$ be a submatrix of $S_4$ with the same number of rows as $S_4$. If $S_4'$ allows o-SIPP, then
\[
S' = \left[\begin{array}{c|c}
S_1 & O \\
\hline
S_3 & S_4'
\end{array}\right]
\]
allows row orthogonality.
\end{proposition}
\begin{proof}
Let $Q$ be a row orthogonal realization of $S$. Then 
\[
Q = 
\left[\begin{array}{c|c}
Q_1 & O \\
\hline
Q_3 & Q_4
\end{array}\right],
\]
where the partition is consistent with that of $S$. Assume that $S_4'$ allows o-SIPP. Then there exists a row orthogonal realization  $Q_4'$ of $S_4'$ with the SIPP. Then $\left[\begin{array}{c|c} O &  Q_4'\end{array}\right]$ is row orthogonal and, by Lemma~\ref{lem: add O}, has the SIPP. By Theorem \ref{thm: gen sup pat}, there exists an $\epsilon > 0$ and a matrix $M_\epsilon'$ such that  $\left[\begin{array}{c|c}\epsilon Q_3 & M_\epsilon'\end{array}\right]$ is row orthogonal and $\sgn(Q_4') = \sgn(M_\epsilon')$. Since $Q$ is row orthogonal, $Q_1 Q_3^T = O$. Therefore,
\[
Q'= \left[\begin{array}{c|c}
Q_1 & O \\
\hline
\epsilon Q_3 & M_\epsilon'
\end{array}\right]
\]
is row orthogonal and $\sgn(Q') = S'$.
\end{proof}

\subsection{Sign patterns requiring o-SIPP}\label{sub:requires}

In this section we present results concerning sign patterns that require o-SIPP. As we shall see, both the number of zero entries and the location of the zero entries in a sign pattern $S$ play an important role in determining whether $S$ requires o-SIPP.

While sign patterns that require o-SIPP have not been previously studied, there are some known results that are closely related to requiring o-SIPP. For example, consider the $n\times n$ lower Hessenberg matrix
\[
H =
\left[\begin{array}{ccccc}
1 & -1 & 0 & \cdots & 0\\
\vdots & \ddots & -2 & \ddots & \vdots\\
\vdots & & \ddots & \ddots & 0\\
\vdots & &  & \ddots & -(n-1)\\
1 & \cdots & \cdots & \cdots & 1
\end{array}\right],
\]
which has orthogonal rows. The proof of Corollary 5.2 in \cite{CS2020} implies that any sign pattern that  has the same zero-nonzero pattern as $H$ and allows orthogonality requires o-SIPP.

The next lemma provides an example of a structural property that guarantees a matrix has the SIPP and is used to establish Corollary \ref{c:lowerHbergSP}, which is  a  slightly more general result than   \cite[Corollary 5.2]{CS2020}. For any integer $k=1-m,\dots,n-1$, the \emph{$k$-th diagonal} of an $m\times n$ matrix $A = [a_{ij}]$ is the list of entries $a_{ij}$ such that $j-i = k$. The $k$-th diagonal terminology is also applied to sign patterns.

\begin{lemma}\label{lem:nz-utriang-has-SIPP}
Let $A = [a_{ij}]\in\R^{m\x n}$ be a wide matrix with full rank. Suppose that there is an integer $k$ such that  $0\le k\le n-1$, each entry of $A$ on  the $r$-th diagonal is nonzero for $ 1-m\le r\le  k$, and   each entry of $A$ on the $r$-th diagonal is zero for $ k< r\le  n-1$. Then $A$ has the SIPP.  
\end{lemma}
\begin{proof}
Note that if $k=n-1$, then $A$ is nowhere zero and hence has the SIPP. Suppose that $k< n-1$. Let $c = \min\{n-k-1,m\}$ and $A_\ell = A[\{1,\ldots,\ell\},\{1\ldots,\ell+k+1\}]$ for $\ell = 1,\ldots, c$. We begin by successively showing that each $A_\ell$ has the SIPP. Since $A_1$ contains a nonzero entry, Lemma \ref{lem:(1,1)-block-has-SIPP} implies $A_1$ has the SIPP. Suppose that $A_i$ has the SIPP for some $i\in\{1,\ldots,c-1\}$. Then Lemma \ref{lem:(1,1)-block-has-SIPP} and Lemma \ref{lem: add O} imply that $A_{i+1}$ has the SIPP.  If $c=m$, then $A=\left[\begin{array}{c|c}A_c & O \end{array}\right]$ has the SIPP by Lemma \ref{lem: add O}. Otherwise, $A = \left[\begin{array}{c} A_{c}\\\hline B \end{array}\right]$, where $B$ is nowhere zero. By Lemma \ref{lem:(1,1)-block-has-SIPP}, $A$ has the SIPP.
\end{proof}

\begin{corollary}\label{c:lowerHbergSP}
Let $S$ be an $m\times n$ wide sign pattern. Suppose that there is an integer $k$ such that  $0\le k\le n-1$, each entry of $S$ on  the $r$-th diagonal is nonzero for $ 1-m\le r\le  k$, and   each entry of $S$ on the $r$-th diagonal is zero for $ k< r\le  n-1$.  If $S$ allows row orthogonality, then $S$ requires o-SIPP.
\end{corollary}

For  sign equivalent matrices $A$ and $B$, $A\in\orth(m,n)$ implies $B\in\orth(m,n)$ and $A$ has the SIPP implies $B$ has the SIPP. Thus the analogous statement with the upper part nonzero is also true.

\begin{corollary}
Let $S$ be a wide $m\x n$ sign pattern. Suppose that there is an integer $k$ such that  $1-m\le k \le 0$, each entry of $S$ on the $r$-th diagonal is nonzero for $ k\le r\le  n-1$, and each entry of $S$ on  the $r$-th diagonal is zero for $ 1-m\le r<  k$. If $S$ allows row orthogonality, then $S$ requires o-SIPP.  
\end{corollary}

In this paper, a \textit{nonzero hollow} matrix (respectively,  sign pattern) is a square matrix (respectively,  sign pattern) with zeros along the main  diagonal and nonzero entries everywhere else. Recall that a \textit{signature matrix} is a diagonal matrix with diagonal entries equal to $\pm1$. Matrices $A$ and $B$ are \textit{signature equivalent} if there exist signature matrices  $D_1$ and $D_2$ such that $D_1 A D_2 = B$. Similarly, sign patterns $S$ and $R$ are signature equivalent if there exist signature matrices  $D_1$ and $D_2$ such that $D_1 C_S D_2 = C_R$. Theorem 5.7 in \cite{CS2020} states that a nonzero hollow matrix $Q\in\orth(n)$ has the SIPP if and only if $Q$ is not signature equivalent to a symmetric hollow matrix. The following corollary is an immediate consequence.

\begin{corollary}\label{cor: hollow req}
Let $S$ be a nonzero  hollow sign pattern that allows orthogonality. If $S$ is not signature equivalent to a symmetric hollow sign pattern, then $S$ requires o-SIPP.  If $S$ is signature equivalent to a symmetric hollow sign pattern, then $S$ does not allow o-SIPP.
\end{corollary}

As an example, consider
\[
S = 
\left[\begin{array}{cccc}
0 & + & + & +\\
+ & 0 & - & +\\
+ & + & 0 & -\\
+ & - & + & 0
\end{array}\right].
\]
It is not difficult to verify that $C_S$ has orthogonal rows, and hence $S$ allows orthogonality. Further, $S$ is not signature equivalent to a symmetric sign pattern. Thus, $S$ requires o-SIPP.

Let $G$ be a graph with vertex set  $V(G)=\{v_1,\dots,v_m\}$ and edge set $E(G)=\{e_1,\dots,e_\ell\}$. The  (vertex-edge) incidence matrix $R_G=[r_{ij}]$ of  $G$ is the $m\x \ell$ matrix that has $r_{ij}=1$ if $v_i\in e_j$ and $r_{ij}=0$ otherwise.  An orientation $\vec G$ of  $G$ is the assignment of a direction to each edge.  That is, the edge $e_j=\{v_i,v_k\}$ is replaced by exactly one of the two arcs $(v_i,v_k)$ or $(v_k,v_i)$; the arc associated with $e_j$ is denoted by $\vec e_j$. The incidence matrix   $R_{\vec G}=[r_{ij}]$ of an orientation $\vec G$ of $G$ is the $m\x \ell$ matrix that has $r_{ij}=-1$ if $\vec e_j=(v_i,v_k)$, $r_{ij}=1$ if $\vec e_j=(v_k,v_i)$, and $r_{ij}=0$ otherwise.

Consider the complete graph $K_m$ and an orientation $\vec K_m$.  For $m\ge 2$, define the $m\x 2 \binom{m}{2}$ matrix $\RKm m=[R_{K_m} | R_{\vec K_m}]$ and its sign pattern  $\SKm m=\sgn(\RKm m)$. The sign pattern $\SKm m$  was shown to have a row orthogonal realization with the SIPP in \cite{Curtis_20}. We now show that $\SKm m$ requires o-SIPP. We note that this sign pattern will be instrumental in Section \ref{sec:asymptotics} for studying random sign patterns that allow row orthogonality.

\begin{theorem}\label{incidence-requires-o-sipp}
For $m\ge 2$ and any   orientation $\vec K_m$ of $K_m$, the sign pattern $\SKm m$ requires o-SIPP.
\end{theorem}
\begin{proof}
For brevity, let $S = \SKm m$. It is not difficult to verify that $C_S$ has orthogonal rows and hence $S$ allows row orthogonality. Let $Q\in \Q(S)$ be row orthogonal and suppose $X=[x_{ij}] \in \sym(m)$ satisfies $(XQ)\circ Q=O$.  By Lemma \ref{lem:identity-circ-Q}, $x_{ii}=0$ for $i=1,\dots,m$. For $i\ne k$, there is a unique edge $e_j=\{ v_i,v_k\}$.  Then applying Lemma \ref{lem:2nz-col} to column $j$ gives that $x_{ik}=0$.  Thus $X = O$ and $Q$ has the SIPP.
\end{proof}

Let $Q\in\orth(m,n)$. It is an immediate consequence of Theorem \ref{t:basic-SIPP} that if $Q$ has a pair of combinatorially orthogonal rows, then $Q$ does not have the SIPP. Similarly, if $m=n$ and $Q$ has a pair of combinatorially orthogonal columns, then $Q$ does not have the SIPP. Thus, when studying sign patterns that require o-SIPP, it is natural to assume that the rows (and in the square case the columns) are not cominbatorially orthogonal. Note that it is possible for a wide sign pattern to have combinatorially orthogonal columns and still require o-SIPP. For example, consider
\[
S =
\left[\begin{array}{cccc}
+ & 0 & + & +\\
0 & + & + & -\\
0 & - & + & -
\end{array}\right]
\quad \text{and} \quad
A =
\left[\begin{array}{cccc}
1 & 0 & 1 & 1\\
0 & 1 & 1 & -1\\
0 & -2  & 1 & -1
\end{array}\right].
\]
Observe that $A$ has orthogonal rows. As we shall see, Corollary \ref{c:allbut3rowsnonzero} implies $S$ requires o-SIPP. 

The remainder of this section investigates how restricting the location and number of zero entries affects whether or not a sign pattern requires o-SIPP. Let $A\in\R^{m\times n}$ be a wide matrix. It is not difficult to verify that if $A$ is nowhere zero, then $A$ has the SIPP. Thus every nowhere zero sign pattern that allows row orthogonality requires o-SIPP. For each additional $0$ entry in $A$, the equation $(XA)\circ A = O$ imposes one fewer linear equation on the entries of $X$. This suggests that the more $0$ entries $A$ has, the larger the solution space to $(XA)\circ A = O$ is going to be, reducing the likelihood that $A$ has the SIPP. In fact, this is the intuition behind the next theorem, which bounds the number of zero entries a matrix with the SIPP can have. 

\begin{theorem}{\rm \cite{CS2020}}\label{thm: bnd 0s}
Let $Q\in\orth(m,n)$ have the SIPP. Then the number of zero entries in $Q$ is at most $nm - \frac12 m(m+1)$.
\end{theorem}

The location of $0$ entries in a sign pattern $S$ also play an important role in determining if $S$ requires o-SIPP.

\begin{theorem}\label{allbut3rowsnonzero}
Let $Q=[q_{ij}] \in \mathcal{O}(m,n)$. Suppose that the zero entries of $Q$ are contained in  at most three rows of $Q$ and that no pair of rows is combinatorially orthogonal. Then, $Q$ has the SIPP.
\end{theorem}
\begin{proof}
Begin by assuming that the zero entries of $Q$ are contained in at most 1 row. Without loss of generality, the zero entries of $Q$ are contained in the first row and the $(1,1)$-entry of $Q$ is nonzero. By Lemma~\ref{lem:(1,1)-block-has-SIPP},  $Q$ has the SIPP.

Now assume that exactly $k$ rows of $Q$ contain a zero, where $k\in\{2,3\}$. Without loss of generality, the first $k$ rows each contain a zero. Let $X \in \sym(m)$ and suppose $(XQ) \circ Q = O$. By Lemma~\ref{lem:zero-rows-of-X}, $X = Y \oplus O$, where $Y\in\sym(k)$.  Observe that $(Y\hat{Q})\circ \hat{Q} = O$, where $\hat{Q}$ is the submatrix of $Q$ formed from the first $k$ rows. Also, note that Lemma \ref{lem:identity-circ-Q} implies $Y\circ I = O$.

First, consider the case $k=2$. Since the rows of $\hat{Q}$ are not combinatorially orthogonal, $\hat{Q}$ has a nowhere zero column. Since $(Y\hat{Q})\circ \hat{Q} = O$, Lemma~\ref{lem:2nz-col} implies that the off-diagonal entries of $Y$ are zero. This, together with $Y\circ I = O$, implies $Y = O$.

Now consider the case $k=3$.  Suppose first that an off-diagonal entry of $Y$, without loss of generality the $(1,2)$-entry, is zero, so 
$Y = 
\left[\begin{array}{ccc}
0 & 0 & y_1 \\
0 & 0 & y_2 \\
y_1 & y_2 & 0 \\
\end{array}\right].
$
Since the rows of $\hat{Q}$ are not combinatorially orthogonal, $\hat Q$ has  a column with  nonzero entries  in the first and third rows. Then $(Y\hat{Q})\circ\hat{Q} = O$ implies $y_1=0$; similarly $y_2=0$.
So suppose that $Y$ is a nonzero hollow matrix. Then by Lemma \ref{lem:2nz-col}, and the preceding argument, no column of $\hat{Q}$ has exactly one zero entry. This, along with the fact that no pair of rows of $\hat{Q}$ are combinatorially orthogonal, implies $\hat{Q}$ has a  nowhere zero column $\bq_j$. Observe that $(Y\hat{Q})\circ\hat{Q} = O$ implies $Y\bq_j = \bzero$. This is impossible since $\rank(Y) = 3$. Thus, $Y = O$ and $Q$ has the SIPP.
\end{proof}

\begin{corollary}\label{c:allbut3rowsnonzero}
Suppose $S$ is an $m\x n$ sign pattern that allows row orthogonality such  that the zero entries of $S$ are contained in at most three rows of $S$ and that no pair of rows are combinatorially orthogonal.
Then, $S$ requires o-SIPP.
\end{corollary}

As the next example illustrates, Corollary \ref{c:allbut3rowsnonzero} cannot be extended to 4 or more rows.

\begin{example}\thlabel{ex:n-zero_n-row}{\rm
Define the $n\times (n+1)$ matrix
\[
A=\left[\begin{array}{cc|cccc|cc}
0 & \sqrt{n-2} & 3-n & 1 & \cdots & 1 & \frac{1}{\sqrt{2}} & \frac{1}{\sqrt{2}}\\
0 & \sqrt{n-2} & 1 & 3-n & \cdots & 1 & \frac{1}{\sqrt{2}} & \frac{1}{\sqrt{2}}\\
\vdots & \vdots & \vdots & \vdots & \ddots & \vdots & \vdots& \vdots\\
0 & \sqrt{n-2} & 1 & 1 & \cdots & 3-n & \frac{1}{\sqrt{2}} & \frac{1}{\sqrt{2}}\\
\hline
0 & \sqrt{n-2} & 1 & 1 & \cdots & 1 & \frac{3-n}{\sqrt{2}} & \frac{3-n}{\sqrt{2}}\\
\hline
\sqrt{n-2} & 0 & 1 & 1 & \cdots & 1 & \frac{1}{\sqrt{2}} & \frac{1}{\sqrt{2}}\\
-\sqrt{n-2} & 0 & 1 & 1 & \cdots & 1 & \frac{1}{\sqrt{2}} & \frac{1}{\sqrt{2}}
\end{array}\right].
\]
Observe that the rows of $A$ are orthogonal. It can be verified that 
\[
X=\bmt
0 & \cdots & 0 & 1 & -1\\
\vdots & \ddots & \vdots & \vdots & \vdots\\
0 & \cdots & 0 & 1 & -1\\
1 & \cdots & 1 & 0 & 0\\
-1 & \cdots & -1 & 0 & 0
\emt\in \R^{n\times n}
\]
satisfies $(XA)\circ A=O$. Since $X=X^T$ and $\rank A=n$,  it follows that $A$ does not have the SIPP.
}
\end{example}

Even if we also prohibit combinatorially orthogonal columns, there are examples of sign patterns with the zero entries restricted to the first four rows that do not require o-SIPP, as seen in the next example, which utilizes a construction method   in \cite{Curtis_20}.

\begin{example}
{\rm
Consider
\[
Q =
\left[
\begin{array}{ccccccc}
 -9 & 9 & 0 & 0 & 3 \sqrt{2} & -6 \sqrt{2} & 6 \sqrt{2} \\
 9 & -9 & 0 & 0 & 3 \sqrt{2} & -6 \sqrt{2} & 6 \sqrt{2} \\
 0 & 0 & -9 & 9 & -6 \sqrt{2} & 3 \sqrt{2} & 6 \sqrt{2} \\
 0 & 0 & 9 & -9 & -6 \sqrt{2} & 3 \sqrt{2} & 6 \sqrt{2} \\
 3 \sqrt{2} & 3 \sqrt{2} & -6 \sqrt{2} & -6 \sqrt{2} & 8 & 8 & 4 \\
 -6 \sqrt{2} & -6 \sqrt{2} & 3 \sqrt{2} & 3 \sqrt{2} & 8 & 8 & 4 \\
 6 \sqrt{2} & 6 \sqrt{2} & 6 \sqrt{2} & 6 \sqrt{2} & 4 & 4 & 2 \\
\end{array}
\right]
\quad \text{and} \quad
X =
\left[\begin{array}{cccc}
0 & 0 & 1 & -1\\
0 & 0 & 1 & -1\\
1 & -1 & 0 & 0\\
1 & -1 & 0 & 0\\
\end{array}\right]\oplus O.
\]
It is readily verified that the rows of $Q$ are orthogonal and $(XQ)\circ Q = O$. Thus, $\sgn(Q)$ allows orthogonality but does not require o-SIPP.
}
\end{example}

By restricting the number of zero entries, we obtain the following result.

\begin{proposition}\label{prop: 4 zeros o-SIPP}
Let $Q\in\orth(m,n)$ have at most four zero entries. Suppose that no pair of rows and no pair of columns columns of $Q$ are combinatorially orthogonal. Then $Q$ has the SIPP.
\end{proposition}
\begin{proof}
Observe that if the zero entries of  $Q = [q_{ij}]$ are contained in at most three rows, then Theorem~\ref{allbut3rowsnonzero} implies $Q$ has the SIPP. So, suppose that each zero entry of $Q$ is contained in a unique row.

Assume  first that $m = n = 4$. If $Q$ has a nowhere zero column, then Proposition \ref{prop: equiv transpose} and Theorem \ref{allbut3rowsnonzero} imply $Q$ has the SIPP. Otherwise, without loss of generality, $Q$ is a nonzero hollow matrix. Theorem 3.2 in \cite{omzd2019} guarantees that $Q$ is not symmetric. Thus, Corollary \ref{cor: hollow req} implies $Q$ has the SIPP. 

Now assume that $n\geq 5$ (and $m\geq 4$). Without loss of generality, $Q$ has the form
\[
Q = 
\left[\begin{array}{c|c}
Q_{1} & Q_{2} \\ 
\hline
Q_{3} & Q_{4}
\end{array}\right],
\]
where all four zero entries of $Q$ are contained in $Q_1$ such that each column and each row of $Q_{1}$ contains a zero entry. Note that $Q_1$ has four rows, $Q_{2}$ is nowhere zero, and $Q_{3}$ and $Q_4$ may be vacuous. We proceed via cases.

We first resolve the case where $Q_1$ has exactly one column, i.e., $Q_1 = \bzero$. Then $Q_2\in\orth(4,n-1)$ is nowhere zero. Thus $Q_2$ has the SIPP and, by Lemma~\ref{lem:(1,1)-block-has-SIPP}, $Q$ has the SIPP.

Suppose that $Q_1$ has at least two columns. Let $X \in \sym(m)$ and suppose $(XQ) \circ Q = O$. By Lemma~\ref{lem:zero-rows-of-X}, $X = Y \oplus O$, where $Y$ is symmetric and has four rows. Further, Lemma \ref{lem:identity-circ-Q} implies $Y\circ I = O$, i.e., $Y$ has the form
\[
Y = 
\left[\begin{array}{cccc}
0 & y_1 & y_2 & y_3\\
y_1 & 0 & y_4 & y_5\\
y_2 & y_4 & 0 & y_6\\
y_3 & y_5 & y_6 & 0
\end{array}\right].
\]

We now consider the case where $Q_1$ has a column $\bc$ with exactly one zero entry. Without loss of generality, the zero appears in the last entry of $\bc$. Then $(Y\bc)\circ \bc = \bzero$ implies $y_1 = y_2 = y_4 = 0$. It now follows from $(YQ_2)\circ Q_2 = O$ that $y_3 = y_5 = y_6 = 0$. Thus, $Y= O$ and $Q$ has the SIPP.

For the final case, assume that $Q_1$ has exactly two columns that contain exactly two zero entries each. Since the first two columns cannot be combinatorialy orthogonal, there must be at least five rows.  Then, without loss of generality, 
\[
\hat{Q}:=
Q[\{1,\ldots,5\},\{1,\ldots,n\}] =
\left[\begin{array}{cc|c}
0 & q_{12} & \bv_1^T\\
0 & q_{22} & \bv_2^T\\
q_{31} & 0 & \bv_3^T\\
q_{41} & 0 & \bv_4^T\\
\hline
q_{51} & q_{52} & \bw^T
\end{array}\right].
\]
From $(YQ_1)\circ Q_1 = O$ it follows that $y_1 = y_6 = 0$.  If $y_i = 0$ for $i\in\{2,3,4,5\}$, then $(YQ_2)\circ Q_2 = O$ implies $Y = O$. Suppose this is not the case, i.e., $y_i\not=0$ for all $i\in\{2,3,4,5\}$. Then $(YQ_2)\circ Q_2 = O$ implies $y_2\bv_1^T+y_4\bv_2^T=\bzero^T$, and thus $\bv_2 = a\bv_1$ for some nonzero value $a$.  Since $Q$ is row orthogonal
\begin{eqnarray*}
0 &=& q_{12}q_{22} + \bv_1^T\bv_2 = q_{12}q_{22} + a(1-q_{12}^2),\\
0 &=& q_{22}q_{52} + \bv_2^T\bw = q_{22}q_{52} + a\bv_1^T\bw, \text{ and}\\
0 &=& q_{12}q_{52} + \bv_1^T\bw.
\end{eqnarray*}
From the last two equations $a = q_{22}/q_{12}$. Substituting this into the first equation implies $q_{22} = 0$, a contradiction. Thus, $Y = O$ and $Q$ has the SIPP.
\end{proof}

It is possible to show that if  $Q\in\orth(m,n)$ has at most five zero entries, no pair of rows and no pair of columns columns of $Q$ are combinatorially orthogonal, then $Q$ has the SIPP. However, with the available tools, the argument is not illuminating and does not warrant the space that would be required. As the next example illustrates, it is possible for a sign pattern with six zero entries to allow row orthogonality, not have combinatorially orthogonal rows or columns, and not require o-SIPP.
\begin{example}
{\rm
Consider the sign pattern
\[
S =
\left[\begin{array}{cccccc}
 0 & + & + & + & + & + \\
 + & 0 & + & - & + & - \\
 + & + & 0 & + & - & - \\
 + & - & + & 0 & - & + \\
 + & + & - & - & 0 & + \\
 + & - & - & + & + & 0
\end{array}\right].
\]
Observe that $C_S$ is a  conference matrix, i.e., $C_S$ is hollow, every off-diagonal entry is $1$ or $-1$, and $C_S^TC_S=(n-1)I$. Hence $S$ allows orthogonality. By Corollary \ref{cor: hollow req}, $S$ does not require o-SIPP.  In fact, it is not difficult to see that $C_S$ does not have the SIPP since the symmetric matrix $X= C_S$ satisfies $(X C_S)\circ C_S = O$. 
}
\end{example}

\subsection{Nowhere zero sign patterns that minimally allow orthogonality}\label{ss:5xn}

In this section we determine the  nowhere zero  sign patterns with at most five rows that minimally allow orthogonality. These and previously known results are summarized in  Table  \ref{t:unique-tar-data} at the end of this section, which lists a representative of each equivalence class of  $m\times n$ nowhere zero sign patterns that minimally allow orthogonality for $m\leq 5$. Recall that a sign pattern $S$ \textit{minimally allows orthogonality} provided $S$ allows row orthogonality and every  sign pattern obtained from $S$ by  deleting  one or more columns  does not allow row orthogonality. 

A complete characterization of nowhere zero sign patterns with at most 4 rows that minimally allow orthogonality  was presented in \cite{Curtis_20}. We summarize these results in the following theorem.

\begin{theorem}  {\rm \cite[Section 5.2]{Curtis_20}}
Let $S$ be an $m\times n$ nowhere zero sign pattern. If $m\leq 3$, then $S$ minimally allows orthogonality if and only if $n=3$ and $S$ is row and column PPO. If $m = 4$, then $S$ minimally allows orthogonality if and only if $n=4$ and $S$ is row and column PPO, or $S$ is sign equivalent to
\[
\left[\begin{array}{ccccc} 
- & - & + & + & + \\
+ & + & - & - & + \\
+ & + & + & + & - \\
+ & + & + & + & +
\end{array} \right].
\]
\end{theorem}

We now determine all $5 \times n$ nowhere zero sign patterns that minimally allow orthogonality.  The next theorem establishes the square case.
\begin{theorem}{\rm \cite[Theorem 7.9]{Curtis_20}}\label{thm:5x5nwz_min}
Let $S$ be a $5\times 5$ nowhere zero sign pattern. Then $S$ allows orthogonality if and only if $S$ is row and column PPO.
\end{theorem}

\begin{lemma}\label{lem:5x4num_zero}
Let $S$ be a $5 \times 4$ nowhere zero sign pattern. Then $S$ is sign equivalent to a sign pattern with at most 5 negative entries.
\end{lemma}
\begin{proof}
By scaling the rows and columns of $S$ we can obtain the sign patterns
\[
S_1 = \left[\begin{array}{c|ccc}
+ & + & + & +\\
\hline
+ & \\
+ & & \multirow{2}{*}{$R$}\\
+ &\\
+ &\\
\end{array}\right]
\quad \text{and} \quad
S_2 = \left[ \begin{array}{c|ccc}
- & + & +&+  \\
\hline
+ & \\ 
+ & & \multirow{2}{*}{$-R$} \\
+ & \\
+ & \\
\end{array} \right]. 
\]
If $R$ contains at most 5 negative entries, then the proof is complete. So, suppose that $R$ has at least 6 negative entries.

First consider the case where $R$ has exactly 6 negative entries.  If a row (or column) of $R$ has 3 negatives, then negating the corresponding row (or column) of $S_1$ reduces the total number of negative entries to at most 5.  Otherwise $R$ has two rows $\br_1$ and $\br_2$, each containing exactly 2 negative entries, and a third row that contains a negative entry; let $j$ denote the column index of this entry. Observe that negating the rows of $S_1$ corresponding to $\br_1$ and $\br_2$ does not change the total number of negative entries. Thus, we can scale the rows of $S_1$ so that column $j$ has 3 negative entries. Negating column $j$ now reduces the total number of negatives to at most 5.

Now suppose that $R$ has at least 7 negative entries. Then $S_2$ has at most 6 negative entries. As before, we can reduce the total number of negative entries to at most 5 if a row (or column) of $-R$ contains at least 3 negative entries, or if $-R$ has two rows that each contain exactly 2 negative entries. Thus, we may assume $-R$ has 1 row with exactly 2 negative entries and 2 columns with exactly 2 negative entries. Without loss of generality $S_2$ is of the form
\[
\left[\begin{array}{c|ccc}
- & + & + & + \\
\hline
+ & - & - & + \\
+ & + & + & - \\
+ & + & -& + \\
+ & - & + & +
\end{array}\right]
\quad \text{or} \quad
\left[\begin{array}{c|ccc}
- & + & + & + \\
\hline
+ & - & - & + \\
+ & + & + & - \\
+ & + & + & - \\
+ & + & - & +
\end{array}\right].
\]
Negate columns 2 and 3, and rows 1 and 3 in the first case. Negate row 2 and column 4 in the second case.\, 
\end{proof}

The next example uses the ideas illustrated in Example \ref{ex:5by6allows}.

\begin{example}\label{obs:5x6 realizations}
{\rm
We can apply Theorem~\ref{thm:approx} to the matrices
\[
A_1={\small \left[\begin{array}{cccccc}
-424 & -297 & 42 & 382 & 424 & 212 \\ 
290 & 48 & -578 & -70 & 247 & 392 \\
126 & 32 & 2 & 536 & -490 & 310 \\ 
466 & 4 & 39 & 404 & 305 & -407 \\ 
49 & 579 & 384 & 12 & 255 & 301
\end{array}\right]} \text{ and }
A_2={\small\left[\begin{array}{cccccc}
-246 & -246 & 369 & 123 & 369 & 123 \\
494 & -254 & 7 & 127 & 7 & 314 \\
174 & 230 & -11 & -421 & 396 & 75 \\
284 & 107 & 414 & 56 & -41 & -392 \\
2 & 477 & 51 & 367 & 69 & 231
\end{array}\right]}
\]
to obtain row orthogonal matrices with the same sign patterns: For $A_1$, the value $\delta = \frac1{268}>.003$ is obtained from row $3$ of $A_1$ and the value $\epsilon={400\over409890583973\sqrt{409890583973}}<0.0007$ is obtained from rows $1$ and $2$.  Thus $\epsilon < \frac 1{5-1}$. Since $\pert_5$ is increasing on its domain, $\pert_5(\epsilon) < \pert_5(0.0007) < 0.003 < \delta$. For $A_2$, the value $\delta = \frac2{477}>.004$ is obtained from row $5$ of $A_2$ and the value $\epsilon={1\over395150\sqrt{118545}}<0.0009$ is obtained from rows $1$ and $2$.  Thus $\epsilon < \frac 1{5-1}$ and $\pert_5(\epsilon) < \pert_5(0.0009) < 0.004 < \delta$.
}
\end{example}

\begin{theorem}
Let $S$ be a $5\times n$ nowhere zero sign pattern. Then $S$ minimally allows orthogonality if and only if $n= 5$ and $S$ is row and column PPO, or $S$ is sign equivalent to 
\[
 S_1={\small \left[\begin{array}{cccccc}
- & - & + & + & + & + \\
+ & + & - & - & + & +\\
+ & + & + & + & - & + \\
+ & + & + & + & + & - \\
+ & + & + & + & + & +
\end{array} \right]},
\ 
S_2={\small\left[\begin{array}{cccccc}
- & - & - & + & + & + \\
+ & + & - & + & + & +\\
+ & + & + & - & - & + \\
+ & + & + & + & + & - \\
+ & + & + & + & + & +
\end{array} \right]},
\ \text{or}
\ 
S_3={\small\left[\begin{array}{cccccc}
- & - & + & + & + & + \\
+ & - & + & + & + & +\\
+ & + & - & - & + & + \\
+ & + & + & + & - & - \\
+ & + & + & + & + & +
\end{array} \right]}.
\]
\end{theorem}
\begin{proof}
Observe that each of the three patterns $S_1, S_2$, and $S_3$  allows row orthogonality by Examples \ref{ex:5by6allows} and \ref{obs:5x6 realizations}.  Removing a column from one of  $S_1, S_2$ or $S_3$ results in a $5\x 5$ sign pattern with a  duplicate column, and such a sign pattern  is not column PPO.  So by Theorem \ref{thm:5x5nwz_min}, removing a column from one of $S_1, S_2$ and $S_3$ results in a sign pattern that does not allow orthogonality. Thus, each of $S_1, S_2$ and $S_3$ minimally allows row orthogonality.

Assume that $S$ minimally allows orthogonality. Without loss of generality the first row and first column of $S$ have all positive entries. Suppose that $S$ has $d$ distinct columns $\bc_1,\ldots,\bc_d$. It is easy to see that $d\geq 4$:  If $S$ had at most 3 distinct columns, then $S$ would have at most $4$  distinct rows, contradicting the fact that $S$ is row PPO.

First consider the case where $S$ has $d = 5$ distinct columns. Since $S$ is row PPO,
\[ 
R = \left[ \begin{array}{c|c|c|c|c}
\bc_1 & \bc_2 & \bc_3 & \bc_4 & \bc_5  \\
\end{array}\right]
\]
is row PPO. Observe that $R$ is column PPO since $\bc_1,\ldots,\bc_5$ are distinct. By Theorem \ref{thm:5x5nwz_min}, $R$ allows orthogonality. Since $S$ minimally allows orthogonality, $S = R$.

Now suppose that $d > 5$ and let  
\[ 
R = \left[ \begin{array}{c|c|c|c|c}
\bc_1 & \bc_2 & \bc_3 & \bc_4 & \bc_5  \\
\end{array}\right].
\]
As before, if $R$ is row PPO, then $R$ allows orthogonality. Since $S$ minimally allows orthogonality, it follows from the preceding argument that $R$ is not row PPO. Without loss of generality
\[ 
R = \left[
\begin{array}{r|rrrr}
+ & + & + & + & +\\
+ & + & + & + & +\\
\hline
+ & \\
+ & & \multicolumn{2}{c}{\hat{R}}\\
+ &\\
\end{array}
\right],
\]
where the columns of $\hat{R}$ are distinct and belong to the set
\[
\left\{
\left[ \begin{array}{ccc}
+ \\
+ \\
-
\end{array} \right], 
\left[ \begin{array}{ccc}
+ \\
- \\
+
\end{array} \right],
\left[ \begin{array}{ccc}
- \\
+ \\
+
\end{array} \right],
\left[ \begin{array}{ccc}
+ \\
- \\
-
\end{array} \right],
\left[ \begin{array}{ccc}
- \\
+ \\
-
\end{array} \right],
\left[ \begin{array}{ccc}
- \\
- \\
+
\end{array} \right],
\left[ \begin{array}{ccc}
- \\
- \\
-
\end{array} \right]
\right\}.
\]
Either $\hat{R}$ contains a column with exactly one negative entry or every column of $\hat{R}$ has at least two negative entries. Observe that in the latter case, negating the last three rows of $R$ results in a column with all positive entries and a column with exactly one negative entry. Thus, we may assume
\[ 
R = \left[ \begin{array}{cc|ccc}
+ & + & + & + & + \\
+ & + & + & + & + \\
\hline
+ & + & + & - & * \\
+ & + & - & * & * \\
+ & - & * & * & *
\end{array}\right].
\]
Since $S$ allows row orthogonality, $\bc_j = (+, -, *,*,*)^{T}$ for some $j\geq 6$. Observe that
\[ 
\left[ \begin{array}{c|c|c|c|c}
\bc_1 & \bc_2 & \bc_3 & \bc_4 & \bc_j  \\
\end{array}\right]
\]
is row and column PPO and hence allows orthogonality by Theorem \ref{thm:5x5nwz_min}. This is a contradiction since $S$ minimally allows orthogonality. Thus, $d\leq 5$.

Finally, we consider the case where $S$ has exactly $d = 4$ distinct columns. Let 
$R = \left[\begin{array}{c|c|c|c}
\bc_1 & \bc_2 & \bc_3 & \bc_4
\end{array}\right]$.
By Lemma \ref{lem:5x4num_zero} we may assume that $R$ has at most $5$ negative entries. Observe that at least $4$ rows of $R$ contain a negative entry since $S$ is row PPO. 

Suppose $R$ has exactly $4$ negative entries. Then $R$ is sign equivalent to
\[
\left[ 
\begin{array}{cccc}
- & + & + & +  \\
+ & - & + & + \\
+ & + & - & + \\
+ & + & + & - \\
+ & + & + & +
\end{array}
\right]
\]
and we assume $R$ has this form. Observe that $S$ can be obtained from $R$ by duplicating some of the columns. By Theorem \ref{thm: big rank 1} we must duplicate at least 2 distinct columns of $R$ to obtain $S$. It follows that, up to sign equivalence, $S$ contains the submatrix
\[
S_1 =
 \left[\begin{array}{cccccc}
- & - & + & + & + & + \\
+ & + & - & - & + & +\\
+ & + & + & + & - & + \\
+ & + & + & + & + & - \\
+ & + & + & + & + & +
\end{array} \right].
\]

Suppose that $R$ has $5$ negatives. Observe that $R$ cannot have exactly 1 negative per row, as this would contradict $S$ being row PPO. Further, at most 1 row of $R$ has 2 negatives, otherwise we have 2 positive rows, which violates row PPO. By these considerations, $R$ is sign equivalent to
\[
\left[\begin{array}{cccc}
- & - & + & +  \\
+ & - & + & +\\
+ & + & - & +\\
+ & + & + & -\\
+ & + & + & +\\
\end{array}\right]
\]
and we assume $R$ has this form.  As before, $S$ can be obtained from $R$ by duplicating at least 2 distinct columns of $R$. Observe that duplicating only columns 1 and 2 of $R$ violates Theorem \ref{thm: big rank 1}. Duplicating columns 1 and 3 is sign equivalent to duplicating columns 1 and 4. Duplicating columns 2 and 3 is sign equivalent to duplicating columns 2 and 4. Thus, up to sign equivalence, $S$ contains one of
\[
S_2 =
\left[\begin{array}{cccccc}
- & - & - & + & + & + \\
+ & + & - & + & + & +\\
+ & + & + & - & - & + \\
+ & + & + & + & + & - \\
+ & + & + & + & + & +
\end{array} \right],
\quad
S_3 = \left[\begin{array}{cccccc}
- & - & + & + & + & + \\
+ & - & + & + & + & +\\
+ & + & - & - & + & + \\
+ & + & + & + & - & - \\
+ & + & + & + & + & +
\end{array} \right]
\]
or
\[
S_4 = \left[\begin{array}{cccccc}
- & - & - & + & + & + \\
+ & - & - & + & + & +\\
+ & + & + & - & - & + \\
+ & + & + & + & + & - \\
+ & + & + & + & + & +
\end{array} \right]
\]
as a submatrix. Observe that $S_4$ is sign equivalent to $S_2$ (negate rows $1,2$ and $5$, and negate columns $4,5$ and $6$; then appropriately permute rows and columns). By Examples \ref{ex:5by6allows} and  \ref{obs:5x6 realizations}, $S_1, S_2$ and $S_3$ allow row orthogonality. Since $S$ minimally allows orthogonality, $S$ is sign equivalent to $S_1, S_2,$ or $S_3$.
\end{proof}

\begin{remark} 
{\rm
Characterizing all $6\times n$ sign patterns that minimally allow orthogonality may require a new approach. However, in doing so, we may learn a great deal about sign patterns that allow row orthogonality. Consider the $6\times 8$ sign pattern
\[
S =
\left[\begin{array}{cccccccc}
+ & + & + & + & + & + & + & + \\ 
+ & + & + & - & - & - & + & + \\ 
+ & + & + & + & + & + & - & + \\ 
+ & + & + & - & - & - & - & - \\ 
+ & + & + & + & + & + & + & - \\ 
+ & + & + & - & - & - & + & -
\end{array}\right].
\]
Deleting any number of columns will contradict Theorem \ref{thm: big rank 1}, so if $S$ allows row orthogonality, then it minimally allows orthogonality. Using the techniques described in Example \ref{ex:5by6allows}, we were unable to find a row orthogonal realization of $S$.  It would be very interesting if this sign pattern does not allow row orthogonality.  It is not too difficult to verify that $S$ satisfies the conditions of Theorem \ref{thm: big rank 1}, so this would unveil a new necessary condition for sign patterns to allow row orthogonality.
}
\end{remark}

\begin{table}[!h]
\begin{center}
\begin{tabular}{l|l} 
Rows & Unique sign patterns (up to sign equivalence) 
\\\hline\\
1 & {\small $\left[\begin{array}{c} +\end{array}\right]$}
\\ & \\\hline\\
2 & {\small $\left[\begin{array}{cc} + & -\\ + & + \end{array}\right]$}
\\ & \\\hline\\
3 & {\small $\left[\begin{array}{ccc} + & - & +\\ + & + & -\\ + & + & + \end{array}\right]$}
\\ & \\
\hline\\
4 &
{\small $\left[\begin{array}{cccc} 
+ & - & + & + \\ 
+ & + & - & + \\ 
+ & + & + & - \\
+ & + & + & + 
\end{array}\right]$,
$\left[\begin{array}{cccc} 
+ & - & - & + \\ 
+ & + & - & + \\ 
+ & + & + & - \\
+ & + & + & + 
\end{array}\right]$,
$\left[\begin{array}{cccc} 
- & + & + & + \\ 
+ & - & + & + \\ 
+ & + & - & + \\
+ & + & + & - 
\end{array}\right]$,
$\left[\begin{array}{ccccc} 
- & + & + & + & + \\ 
+ & - & + & - & + \\ 
+ & + & - & + & - \\
+ & + & + & + & +
\end{array}\right]$}
\\ & \\
\hline\\
 &
{\small $\left[\begin{array}{rrrrr} 
+ & - & - & + & +\\
+ & + & - & - & +\\
+ & + & + & - & -\\
+ & + & + & + & -\\
+ & + & + & + & +
\end{array}\right]$,
$\left[\begin{array}{rrrrr} 
- & - & + & + & +\\
+ & - & - & + & +\\
+ & + & - & - & +\\
+ & + & + & + & -\\
+ & + & + & + & +
\end{array}\right]$,
$\left[\begin{array}{rrrrr} 
- & - & + & + & +\\
+ & - & - & + & +\\
+ & + & - & + & +\\
+ & + & + & - & +\\
+ & + & + & + & -
\end{array}\right]$},
\\ & \\
\multirow{7}{1em}{5}&
{\small $\left[\begin{array}{rrrrr} 
- & - & + & + & +\\
+ & - & - & + & +\\
+ & + & + & - & -\\
+ & + & + & + & -\\
+ & + & + & + & +
\end{array}\right]$,
$\left[\begin{array}{rrrrr} 
- & - & + & + & +\\
+ & - & - & + & +\\
+ & + & + & - & +\\
+ & + & + & + & -\\
+ & + & + & + & +
\end{array}\right]$,
$\left[\begin{array}{rrrrr} 
+ & - & - & + & +\\
+ & + & - & + & +\\
+ & + & + & - & -\\
+ & + & + & + & -\\
+ & + & + & + & +
\end{array}\right]$},
\\ & \\
&
{\small $\left[\begin{array}{rrrrr} 
- & - & + & + & +\\
+ & - & + & + & +\\
+ & + & - & + & +\\
+ & + & + & - & +\\
+ & + & + & + & -
\end{array}\right]$,
$\left[\begin{array}{rrrrr} 
- & + & + & + & +\\
+ & - & + & + & +\\
+ & + & - & + & +\\
+ & + & + & - & +\\
+ & + & + & + & -
\end{array}\right]$,
$\left[\begin{array}{rrrrr} 
+ & - & + & + & +\\
+ & + & - & + & +\\
+ & + & + & - & +\\
+ & + & + & + & -\\
+ & + & + & + & +
\end{array}\right]$,}
\\ & \\
&
{\small $\left[\begin{array}{cccccc}
- & - & + & + & + & + \\
+ & + & - & - & + & +\\
+ & + & + & + & - & + \\
+ & + & + & + & + & - \\
+ & + & + & + & + & +
\end{array} \right]$,
\ 
$\left[\begin{array}{cccccc}
- & - & - & + & + & + \\
+ & + & - & + & + & +\\
+ & + & + & - & - & + \\
+ & + & + & + & + & - \\
+ & + & + & + & + & +
\end{array} \right]$,
$\left[\begin{array}{cccccc}
- & - & + & + & + & + \\
+ & - & + & + & + & +\\
+ & + & - & - & + & + \\
+ & + & + & + & - & - \\
+ & + & + & + & + & +
\end{array} \right]$}
\end{tabular}
\end{center}
\caption{One representative of each sign-equivalence class of $m\x n$ nowhere zero sign patterns that minimally allow orthogonality for $m\le 5$}
\label{t:unique-tar-data}
\end{table}
\newpage

%
%
\section{Likelihood a random sign pattern allows row orthogonality}\label{sec:asymptotics}

The question of finding the probability that $m$ vectors sampled from $\{\pm 1\}^n$ are linearly independent has attracted recent attention in the literature. This problem can equivalently be stated as asking for the probability that a random matrix in $\{\pm 1\}^{m \times n}$ has rank $m$ (the literature is most interested in the case when $m=n$). In particular, Tikhomirov~\cite{Tikhomirov2018} answered this question in a strong form by showing that this probability is bounded below by $1-\bigl({1\over 2}+o(1)\bigr)^m$ whenever $n\geq m$; in particular, when $n\geq m$ the probability tends toward $1$ as $m$ tends toward $\infty$.

In this section, we consider the adjacent problem of determining the threshold $t(m)$ such that a random matrix in $\{+,-\}^{m\times n}$ with $n\geq t(m)$ allows row orthogonality with probability tending toward $1$ as $m$ tends toward $\infty$.

Let $f(n)$ and $g(n)$ be functions from the non-negative integers to the reals. Then $f(n) = o(g(n))$ if $\lim_{n\to\infty}{f(n)}/{g(n)} = 0$, and $f(n) = \omega(g(n))$ if $g(n) = o(f(n))$. An event $E=E(n)$ happens \textit{with high probability} as $n\to\infty$ if $\Pr[E]=1-o(1)$. The \emph{union bound} is the fact that the probability  that at least one of a set of the events happens  is at most the sum of the probabilities of the events.

For a probability distribution $\mu$ on a set $\Omega$, we write $x\sim\mu$ to mean that $x$ is distributed according to $\mu$. If $\Omega$ is a finite set, then we write $x\sim\Omega$ to mean that $x$ is chosen uniformly from $\Omega$. We write $x_1,\dots,x_n\sim\mu$ to indicate that $x_1,\dots,x_n$ are distributed according to $\mu$ and are mutually independent from one another. For a positive integer $n$, $\mu^n$ denotes the product distribution on $\Omega^n$ where each entry is drawn independently from $\mu$. Similarly, for positive integers $m$ and $n$, $\mu^{m\times n}$ denotes the product distribution on $\Omega^{m\times n}$ where each entry is drawn independently from $\mu$.  For an index set $\alpha$, let $\Omega^\alpha$ denote  the set of vectors with entries in $\Omega$ indexed by $\alpha$. We write $\mu^\alpha$ to mean the product distribution on $\Omega^\alpha$ where each entry is drawn independently from $\mu$. 

We will need two forms of the Chernoff bound, which we state here.

\begin{theorem}\label{t:Chernoff1}
{\rm \cite[Corollary A.1.2]{probmethod}}
Let $X_i$, $1\leq i \leq n$, be mutually independent random variables with
$
\Pr[X_i = 1] = \Pr[X_i = -1] = \frac12
$
and 
$
X = X_1 + \cdots + X_n.
$
Let $a > 0$. Then
\[
\Pr[|X| > a] < 2e^{-a^2/2n}.
\]
\end{theorem}

\begin{theorem}\label{t:Chernoff2}
{\rm \cite[Remark 9.2]{Tijms}}
Suppose $X_1,\ldots,X_n$ are independent random variables taking values from the set $\{0,1\}$. Let $X = X_1 +\ldots + X_n$. Then for any $\delta \geq 0$
\[
\Pr[X\leq (1-\delta)\mathbb{E}[X]] \leq \exp(-\delta^2 \mathbb{E}[X]/2).
\]
\end{theorem}

For a fixed $0\leq p\leq 1/2$, let $\mu_p$ denote the distribution on $\{0,\pm 1\}$ where $\mu_p(1)=\mu_p(-1)=p$ and $\mu_p(0)=1-2p$. The main result of this section, Theorem~\ref{thm:main prob allows}, implies that for any fixed $0<p\leq 1/2$, there is a constant $C=C(p)$ such that if $A\sim{ \mu_p}^{m\times n}$ where $n \geq m^2 +C m\log m$, then $\sgn(A)$ allows row orthogonality with high probability as $m\to\infty$. Before proving Theorem~\ref{thm:main prob allows}, we use Theorem~\ref{c:close-enough} to obtain a slightly weaker result for nowhere zero sign patterns. We include this result since the proof is relatively short and highlights a substantially different approach from Theorem~\ref{thm:main prob allows}.

\begin{theorem}\label{thm:realization}
If $A\sim\{\pm 1\}^{m\times n}$ with $n\geq 17m^2\log m$, then $\sgn(A)$ allows row orthogonality with high probability as $m\to\infty$.
\end{theorem}
\begin{proof}
Let $\bx_i$ denote the $i$th row of $A$, so $\bx_1,\dots,\bx_m\sim\{\pm 1\}^n$. Observe that $\norm{\bx_i}_2 = \sqrt n$ and that $\delta(\bx_i)=1$ for each $i\in[m]$.
Set
\[
\epsilon=\sqrt{{\frac{17}{4}\log m\over n}},
\]
and observe that $0\leq\epsilon\leq {1\over 2m}$ since $n\geq 17m^2\log m$ and so
\[
r_m(\epsilon) = \sqrt{{1+\epsilon\over (1-(m-2)\epsilon)(1-(m-1)\epsilon)}}-1 \leq \sqrt{{1+{1\over 2m}\over \bigl(1-{m-2\over 2m}\bigr)\bigl(1-{m-1\over 2m}\bigr)}}-1
=2\sqrt{{m(m+1/2)\over(m+2)(m+1)}}-1<1.
\]
Thus, if $\abs{\langle\bx_i,\bx_j\rangle}\leq\epsilon n$ for all $i\neq j\in[m]$, then we may apply Theorem~\ref{c:close-enough} to locate a set of orthogonal vectors $\squig\bx_1,\dots,\squig\bx_m$ such that $\sgn(\squig\bx_i)=\sgn(\bx_i)$. Thus, in order to conclude the proof, it suffices to show that $\abs{\langle\bx_i,\bx_j\rangle}\leq\epsilon n$ for all $i\neq j\in[m]$ with high probability as $m\to\infty$.

Since $\bx_1,\dots,\bx_m\sim\{\pm 1\}^n$ are independent, we may apply the Chernoff bound in Theorem \ref{t:Chernoff1} to bound
\[
\Pr\bigl[\abs{\langle \bx_i,\bx_j\rangle}>\epsilon n\bigr]<2e^{-\epsilon^2n/2}
\]
for any $i\neq j\in\{1,\dots,m\}$. By the  union bound,
\begin{align*}
\Pr\bigl[\abs{\langle \bx_i,\bx_j\rangle}\leq \epsilon n,\ \forall i\neq j\in[m]\bigr] 
&\geq 1 - \sum_{i<j\in[m]}\Pr\bigl[\abs{\langle\bx_i,\bx_j\rangle}>\epsilon n\bigr]\geq 1-\binom{m}{2} 2e^{-\epsilon^2 n/2}\\
& \geq 1-m^2e^{-\epsilon^2n/2} = 1-m^{-1/8} = 1 - o(1).
\end{align*}
\end{proof}

We now show how to improve  Theorem~\ref{thm:realization}  by using the SIPP. Recall that Theorem~\ref{incidence-requires-o-sipp} states that the sign pattern $\SKm m$ requires o-SIPP.  
We say that a pair of negative $4$-cycles are \textit{column-disjoint} if the column indices of the negative $4$-cycles are all distinct.  Observe that any sign pattern that has a collection of column-disjoint negative $4$-cycles between every pair of rows is sign equivalent to a superpattern of $\left[\begin{array}{c|c}  \SKm m & O\end{array}\right]$.  So by Theorem \ref{t:basic-SIPP} and Theorem~\ref{thm: SIPP add zero block}, we have the following  observation.

\begin{observation}
\label{obs: col-disj coll}
Let $S$ be an $m\times n$ sign pattern. If $S$ has a collection of column-disjoint negative 4-cycles between every pair of rows, then $S$ allows row orthogonality. 
\end{observation}

In the following proofs, we must condition on the outcome of a stochastic process.
For those readers unfamiliar with these ideas, we recommend consulting~\cite[Chapter 9]{williams_prob}. 

\begin{lemma}\label{lem:singlerow}
Fix any $0<p\leq 1/2$. If $\bx_1,\dots,\bx_m\sim\mu_p^n$, then the probability that we can find distinct integers $i_1,j_1,i_2,j_2,\dots,i_m,j_m\in[n]$ such that $(\bx_k)_{i_k}=1$ and $(\bx_k)_{j_k}=-1$ for all $k\in[m]$ is at least
\[
1-{1\over p}(1-p)^{n-2m+1}
\]
\end{lemma}
\begin{proof}
We employ the following greedy algorithm  to find the required set $W_m=\{i_1,j_1,\dots,i_m,j_m\}$ of indices:

Initialize $U_0=[n]$  and $W_0=\emptyset$. At time $t=1, \dots, m$, do the following:
\begin{enumerate}
\item 
Reveal $\bx_t$.
\item 
Attempt to locate some $i_t,j_t\in U_{t-1}$ for which $(\bx_t)_{i_t}=1$ and $(\bx_t)_{j_t}=-1$. If such $i_t,j_t$ are found, then set $U_t = U_{t-1}\setminus\{i_t,j_t\}$  and $W_t = W_{t-1}\cup\{i_t,j_t\}$. If no such $i_t,j_t$ exist, then exit with failure.
\end{enumerate}
If the above algorithm succeeds, then we have located the desired $  W_m$.

Let $\tau$ be the  round on which the algorithm fails, setting $\tau=m+1$ if the algorithm succeeds. In order to complete the proof, we show that
\[
\Pr[\tau\leq m]\leq{1\over p}(1-p)^{n-2m+1}.
\]

Fix any $t\in[m]$ and consider conditioning on the event $\{\tau\geq t\}$. Since $\tau\geq t$ if and only if the algorithm has succeeded locating the set $U_{t-1}$, we may condition on such an outcome. Now, conditioned on the algorithm locating $U_{t-1}$, we observe that $\tau=t$ if and only if  $\bx_t[U_{t-1}]$  is nonnegative or nonpositive. Furthermore, before the $t$th loop, no information is known about the vector $\bx_t$ and so $\bx_t[U_{t-1}] \sim \mu_p{}^{U_{t-1}}$. We may therefore bound
\[
\Pr[\tau=t\ \vert\ U_{t-1}] = \Pr\bigl[\bx_t[U_{t-1}]\in\{0,1\}^{U_{t-1}}\cup\{0,-1\}^{U_{t-1}}\ \bigm\vert\ U_{t-1}\bigr]\leq 2(1-p)^{\abs{U_{t-1}}}=2(1-p)^{n-2(t-1)}.
\]
Since this bound is independent of $U_{t-1}$, we may bound
\[
\Pr[\tau=t]\leq\Pr[\tau=t\ \vert\ \tau\geq t]\leq 2(1-p)^{n-2(t-1)}.
\]
We therefore conclude that
\bea
\Pr[\tau\leq m] &= &\sum_{t=1}^m \Pr[\tau=m] \leq \sum_{t=1}^m 2(1-p)^{n-2(t-1)}  =2\sum_{k=0}^{m-1} (1-p)^{n-2m+2+2k}\\ 
&=&2{(1-p)^{n-2m+2}-(1-p)^{n+2}\over (2-p)p}\le {2\over (2-p)p}(1-p)^{n-2m+2}\leq {1\over p}(1-p)^{n-2m+1},
\eea
where the final inequality follows from the fact that ${1-p\over 2-p}\leq {1\over 2}$.
\end{proof}

We now use the above lemma to locate collections of column-disjoint negative $4$-cycles.

\begin{lemma}\label{lem:singlerow 2}
Fix any $\bx \in\{\pm 1\}^n$ and any $0 < p \leq 1/2$. Assume $A\sim\mu_p^{m\times n}$ and set
$
B=\left[\begin{array}{c}
\bx^T \\
\hline
A
\end{array}\right].
$
Then $\sgn(B)$ contains column-disjoint negative $4$-cycles between its first row and all other rows with probability at least
\[
1-{1\over p}(1-p)^{n-2m+1}.
\]
\end{lemma}
\begin{proof}
Let $D$ be the diagonal matrix whose $i$th diagonal entry is the $i$th entry of $\bx$. Observe that $\sgn(B)$ is sign equivalent to
$\sgn( BD)=\sgn\left(\left[\begin{array}{ccc}
1 & \cdots & 1 \\
\hline
& AD &
\end{array}\right]\right)$. 
Since $\mu_p(1) = \mu_p(-1)$ and $\bx\in\{\pm 1\}^n$, $AD\sim\mu_p^{m\times n}$. Thus $\sgn(B)$ contains column-disjoint negative $4$-cycles between its first row and all other rows if and only if we can locate distinct $i_1,\dots,i_m,j_1,\dots,j_m\in[n]$ so that $(\bw_k)_{i_k}=1$ and $(\bw_k)_{j_k}=-1$. As such, the conclusion follows from Lemma~\ref{lem:singlerow}.
\end{proof}

\begin{lemma}\label{lem:col-disj 4cycles}
Fix a number $0<p\leq 1/2$. Assume $A\sim\mu_p^{m\times n}$, where $n\geq m^2+mr+{2m\over p}$ for some  $r\in\{0,\dots,m\}$. Then the probability that $\sgn(A)$ contains a collection of column-disjoint negative $4$-cycles between every pair of rows is bounded below by
\[
1-me^{-m/8}-{m\over p}(1-p)^r.
\]
\end{lemma}
\begin{proof}
In order to locate a collection of negative $4$-cycles between every pair of rows of $\sgn(A)$, we employ the following greedy algorithm: 

Suppose that the rows of $A$ are $\bx_1,\dots,\bx_m$. Initialize $U_0=[n]$. At time $t=1,\dots,m-1$ do the following:
\begin{enumerate}
\item 
Reveal $\bx_t[U_{t-1}]$.
\item\label{step:support} 
Find some $W_t\subseteq \supp( \bx_i[U_{t-1}])$ with $\abs{W_t}=2(m-t)+r$ and set $U_t = U_{t-1}\setminus W_t$. If no such $W_t$ exists, then fail.
\item 
Reveal $A[\{t+1,\dots,m\}, W_t]$
\item\label{step:cycles} 
Locate column-disjoint negative $4$-cycles in $\sgn(A)$ between row $t$ and all rows $k>t$, all of whose columns reside within $W_t$. If such negative $4$-cycles cannot be found, then fail.
\end{enumerate}
If the above algorithm succeeds, then $\sgn(A)$ contains a collection of column-disjoint negative $4$-cycles between every pair of rows.

Let $\tau$ denote the first round on which the algorithm fails, setting $\tau=m$ if the algorithm succeeds. In order to prove the claim, we argue that
\[
\Pr[\tau\leq m-1]\leq me^{-m/8}+{m\over p}(1-p)^r.
\]

Fix any $t\in[m-1]$. Let $\mathcal S_t$ denote the event that the algorithm fails at step~\ref{step:support} on the $t$th loop, and let $\mathcal F_t$ denote the event that the algorithm fails at step~\ref{step:cycles} on the $t$th loop. Certainly $\{\tau=t\}=\mathcal S_t\cup\mathcal F_t$. We begin by bounding the probability of $\mathcal S_t$.

Consider conditioning on the event $\{\tau\geq t\}$. Of course, if $\tau\geq t$, then the algorithm has succeeded in locating the set $U_{t-1}$. Furthermore, conditioned on $\{\tau\geq t\}$ and $U_{t-1}$, observe that prior to the $t$th loop of the algorithm, no entries within $A[\{t,\dots,m\}, U_{t-1}]$ have been revealed; therefore $A[\{t,\dots,m\}, U_{t-1}]\sim\mu_p^{\{t,\dots,m\}\times U_{t-1}}$. In particular, $\bx_t[U_{t-1}]\sim\mu_p^{U_{t-1}}$. Now, $W_t$ cannot be located if and only if $\abs{\supp(\bx_t[U_{t-1}])}<2(m-t)+r$. Additionally,
\bea
\abs{U_{t-1}} &=&n-\sum_{j=1}^{t-1}\bigl(2(m-j)+r\bigr)\geq m^2+rm+{2m\over p}-\sum_{j=1}^{t-1}\bigl(2(m-j)+r\bigr)\\
&=&{2m\over p}+m^2+rm-(t-1)(2m-t)-(t-1)r\geq {2m\over p}.
\eea
We can therefore fix a subset $U\subseteq U_{t-1}$ with $\abs U = { \lf{2m\over p}\rf}$.  From the earlier observation, we know that $\bx_t[U]\sim\mu_p^U$ and so we may bound
\bea
\Pr[\mathcal S_t\ \vert\ \{\tau\geq t\}, U_{t-1}] &=&\Pr\bigl[\abs{\supp(\bx_t[U_{t-1}])}<2(m-t)+r\ \bigm\vert\ \{\tau\geq t\}, U_{t-1}\bigr]\\
&\leq&\Pr\bigl[\abs{\supp(\bx_t[U])}<2(m-t)+r\ \bigm\vert\ \{\tau\geq t\}, U_{t-1}\bigr]\\
&=&\Pr_{\bx\sim\mu_p^{\lfloor 2m/p\rfloor}}\bigl[\abs{\supp\bx}< 2(m-t)+r\bigr] \\
&\leq& \Pr_{\bx\sim\mu_p^{\lfloor 2m/p\rfloor}}\bigl[\abs{\supp\bx}< 3m-1\bigr],
\eea
where the final inequality follows from the fact that $t\geq 1$ and $r\leq m$.

Next, we note that if $\bx\sim\mu_p^{\lfloor 2m/p\rfloor}$, then ${ \mathbb E[\abs{\supp\bx}]}=2p\lf{2m\over p}\rf\geq 4m-2p\geq 4m-1$ and so 
\bea
\Pr_{\bx\sim\mu_p^{\lf 2m/p\rf}}\bigl[\abs{\supp\bx}<3m-1\bigr]&\leq&\Pr_{\bx\sim\mu_p^{\lf 2m/p\rf}}\bigl[\abs{\supp\bx}<{ \mathbb E[\abs{\supp\bx}]}-m\bigr]\\ &=&\Pr_{\bx\sim\mu_p^{\lf 2m/p\rf}}\lb\abs{\supp\bx}<\lp1-\frac{m}{\mathbb E[\abs{\supp\bx}]}\rp\mathbb E[\abs{\supp\bx}]\rb \\
&\leq& \exp\biggl(-\frac{m^2}{2 \mathbb E[|\supp\bx|]}\biggr)
= \exp\biggl(-{m^2\over 4p\lf{2m\over p}\rf}\biggr)\leq
{e^{-m/8}},
\eea
where the second inequality follows from the Chernoff bound {in Theorem \ref{t:Chernoff2}}.  Since this bound is independent of $U_{t-1}$, we have argued that
\begin{equation}\label{eqn:st}
\Pr[\mathcal S_t\ \vert\ \tau\geq t]\leq e^{-m/8}.
\end{equation}

Next, we bound the probability of $\mathcal F_t$. In order for $\mathcal F_t$ to hold, it must be the case that $\tau\geq t$ and that $\mathcal S_t$ does not hold; in particular, the algorithm must have succeeded in locating the set $W_t$. By construction, just after locating $W_t$, no entries within $A[\{t+1,\dots,m\}, W_t]$ have been revealed; therefore $A[\{t+1,\dots,m\}, W_t]\sim\mu_p^{\{t+1,\dots,m\}\times W_t}$. {Since $\Pr[\mathcal F_t\ \vert\ W_t]$ is equal to the probability of not finding a collection of column disjoint negative 4-cycles between the first row of $A[\{t+1,\dots,m\}, W_t]$ and the remaining rows,} we may appeal to Lemma~\ref{lem:singlerow 2} to bound
\[
\Pr[\mathcal F_t\ \vert\ W_t] \leq {1\over p}(1-p)^{\abs{W_t}-2(m-t)+1}\leq{1\over p}(1-p)^r.
\]
Since this bound is independent of $W_t$, we have shown that
\begin{equation}\label{eqn:ft}
\Pr[\mathcal F_t\ \vert\ \{\tau\geq t\}, \overline{\mathcal S_t}]\leq {1\over p}(1-p)^r,
\end{equation}
where $\overline{\mathcal S_t}$ denotes the event that $\mathcal S_t$ does not occur.

Combining \eqref{eqn:st} and \eqref{eqn:ft} we have shown that
\bea
\Pr[\tau=t] &\leq& \Pr[\tau=t\ \vert\ \tau\geq t]=\Pr[\mathcal S_t\ \vert\ \tau\geq t]+\Pr[\mathcal F_t\ \vert\ \tau\geq t]\\
&\leq& \Pr[\mathcal S_t\ \vert\ \tau\geq t]+\Pr[\mathcal F_t\ \vert\ \{\tau \geq t\}, \overline{\mathcal S_t}]\le  e^{-m/8}+{1\over p}(1-p)^r,
\eea
where the first equality holds since $\mathcal S_t$ and $\mathcal F_t$ partition $\{\tau = t\}$. \vspace{-5pt}

Using this inequality, we finally bound\vspace{-5pt}
\[
\Pr[\tau\leq m-1] = \sum_{t=1}^{m-1}\Pr[\tau=t]\leq\sum_{t=1}^{m-1}\lp e^{-m/8}+{1\over p}(1-p)^r\rp\leq me^{-m/8}+{m\over p}(1-p)^r,\vspace{-5pt}
\]
as needed.
\end{proof}

\begin{theorem}\label{thm:main prob allows}
For any fixed $0<p\leq 1/2$, if $A\sim\mu_p^{m\times n}$ and 
\[
n\geq m^2+m\log_{1/(1-p)} m+\omega(m),
\]
then $\sgn(A)$ allows row orthogonality with high probability as $m\to\infty$.
\end{theorem}
\begin{proof}
Suppose that
\[
n\geq m^2+m\log_{1/(1-p)}m+f(m),
\]
where $f(m)=\omega(m)$. Without loss of generality, we may additionally suppose that $f(m)=o(m^2)$. Set $r=\log_{1/(1-p)} m+\frac{f(m)}m-{2\over p}$ which is certainly bounded above by $m$ for all sufficiently large $m$ since $f(m)=o(m^2)$.  Furthermore, by decreasing $f(m)$  by  some amount no more than  $m$, we may ensure that $r$ is an integer,  the lower bound on $n$ remains true,  $f(m)=\omega(m)$, and  $f(m)=o(m^2)$. Now, since $f(m)=\omega(m)$ and $0<p\leq 1/2$ is fixed, we have that
\[
n\geq m^2+m\log_{1/(1-p)}m+f(m)\geq m^2+mr+{2m\over p}
\]
for all sufficiently large $m$. Thus, we may apply Lemma~\ref{lem:col-disj 4cycles} to learn that $\sgn(A)$ contains a collection of column-disjoint negative $4$-cycles between every pair of rows (and hence has a row orthogonal realization) with probability at least
\[
1-me^{-m/8}-{m\over p}(1-p)^r=1-me^{-m/8}-{m\over p}(1-p)^{-\log_{(1-p)}m}(1-p)^{{f(m)\over m}-{2\over p}}=1-me^{-m/8}-{1\over p}(1-p)^{{f(m)\over m}-{2\over p}},
\]
which tends to $1$ as $m\to\infty$ since $f(m)=\omega(m)$ and $0<p\leq 1/2$.
\end{proof}

We suspect that Theorem~\ref{thm:main prob allows} is not best possible.

\begin{question} Determine the threshold $t(m)$ such that if $S\sim\{+,- \}^{m\times n}$ with $n\geq t(m)$, then $S$ has a row orthogonal realization with high probability as $m\to\infty$.
\end{question}

Theorem~\ref{thm:main prob allows} implies that $t(m)\leq m^2+m\log_2m+\omega(m)$. Observe that $t(m)\geq m$ and it is possible that this is the correct answer. As shown in the next theorem, the best known obstruction (see Theorem \ref{thm: big rank 1}) does not block $t(m)=m$.

\begin{theorem}
Let $X\sim \{\pm 1\}^{m\times m}$. Then with high probability as ${m}\to\infty$ the matrix $X$ does not contain an $r\times s$ submatrix $Y$ such that $r + s = {m} + 2$ and $\rank Y = 1$.
\end{theorem}
\begin{proof}
Let $\Omega$ denote the set of pairs $(\bx,\by)$, where $\bx\in \{\pm1\}^r$, $\by\in\{\pm1\}^s$ and the first entry of $\bx$ is $1$. Observe that the map $(\bx,\by)\mapsto \bx\by^T$ is a bijection between $\Omega$ and the set of rank 1 matrices in $\{\pm1\}^{r\times s}$. Thus the probability that $Y\sim\{\pm1\}^{r\times s}$ has rank 1 is precisely $2^{-(r-1)(s-1)}$.

The number of $r\times s$ submatrices of $X$ is $\binom{m}{r}\binom{m}{s}$. By the union bound, the probability that $X$ contains an $r\times s$ submatrix $Y$ such that $r + s = {m} + 2$ and $\rank Y = 1$ is at most
\[
\sum_{r=2}^m \binom{m}{r}\binom{m}{m+2-r}2^{-(r-1)(m+1-r)}
=
\sum_{k=1}^{m-1} \binom{m}{k+1}\binom{m}{m+1-k}2^{-k(m-k)}.
\]
We show that this sum tends toward $0$ as $m \to \infty$ by showing that
\[
\binom{m}{k+1}\binom{m}{m+1-k} < \frac{2^{k(m-k)}}{m^2}
\]
for all $1 \leq k \leq m - 1$, provided $m$ is sufficiently large. If $k \leq 2$ or $k\geq m-2$, then for $m$ sufficiently large
\[
\binom{m}{k+1}\binom{m}{m+1-k} \leq m^4 < \frac{2^{k(m-k)}}{m^2}.
\]
For $3 \leq k\leq m-3$, we have
\[
\binom{m}{k+1}\binom{m}{m+1-k} \leq 2^{2m} < \frac{2^{3(m-3)}}{m^2} \leq \frac{2^{k(m-k)}}{m^2}.
\]
\end{proof}

\section*{Acknowledgements}
The research of Z.~Brennan, C.~Cox, B.~Curtis, E.~Gomez-Leos, and C.~Thompson  was partially supported by NSF grant  1839918 and the authors thank  the National Science Foundation.


\end{document}